\documentclass[10pt]{article}
\usepackage[utf8]{inputenc}
\usepackage{authblk}
\usepackage{amsfonts}
\usepackage{graphicx}
\usepackage{bm}
\usepackage{tikz}
\usepackage{enumitem}
\usepackage[authoryear]{natbib}
\usepackage{multirow}
\usepackage{amsmath}
\usepackage{amsthm}
\usepackage{amssymb}
\usepackage{hyperref}
\hypersetup{
    colorlinks=true,
    citecolor=blue,
    urlcolor=blue,
    linkcolor=blue
}
\usepackage{rotating}
\usepackage{float}
\usepackage{dirtytalk}
\usepackage{lscape}
\usepackage{appendix}
\newtheorem{theorem}{Theorem}[section]
\newtheorem{corollary}{Corollary}[theorem]
\newtheorem{lemma}[theorem]{Lemma}
\theoremstyle{remark}
\newtheorem{remark}{\textbf{Remark}}[section]
\newtheorem{example}{\textbf{Example}}[section]
\usepackage[a4paper, total={6in, 9in}]{geometry}

\usepackage{subcaption}

\newcommand{\sr}[1]{\textcolor{black}{#1}}
\DeclareMathOperator*{\argmin}{arg\,min}

\title{Asymptotic breakdown point analysis of the minimum density power divergence estimator under independent non-homogeneous setups}

\author[1]{\Large{Suryasis Jana}}
\author[2]{\Large{Subhrajyoty Roy}}
\author[1]{\Large{Ayanendranath Basu}}
\author[1]{\Large{Abhik Ghosh\footnote{Corresponding author: abhik.ghosh@isical.ac.in}}}
\affil[1]{Indian Statistical Institute, Kolkata, India}
\affil[2]{Washington University at St. Louis, MO, USA}

\date{}

\begin{document}
\maketitle
% \large

\begin{abstract}
The minimum density power divergence estimator (MDPDE) has gained significant attention in the literature of robust inference due to its strong robustness properties and high asymptotic efficiency; it is relatively easy to compute and can be interpreted as a generalization of the classical maximum likelihood estimator. It has been successfully applied in various setups, including the case of independent and non-homogeneous (INH) observations that cover both classification and regression-type problems with a fixed design. While the local robustness of this estimator has been theoretically validated through the bounded influence function, no general result is known about the global reliability or the breakdown behavior of this estimator under the INH setup, except for the specific case of location-type models. In this paper, we extend the notion of asymptotic breakdown point from the case of independent and identically distributed data to the INH setup and derive a theoretical lower bound for the asymptotic breakdown point of the MDPDE, under some easily verifiable assumptions. These results are further illustrated with applications to some fixed design regression models and corroborated through extensive simulation studies.

\end{abstract} 

% \tableofcontents
\section{Introduction}
In various statistical inference problems, the minimum divergence estimation technique has been observed to be an extremely useful method that is based on the minimization of a suitable statistical distance between the assumed model distribution/density and the true underlying distribution/density of the observed data. We will, in particular, focus on the density-based minimum divergence method. The associated minimum divergence estimator (MDE) is defined to be the value of the parameter that minimizes the chosen divergence between the model density and an empirical estimate of the true data density over the parameter space. The MDE constitutes a super-class containing many popular estimators, such as the (non-robust) maximum likelihood estimator (MLE), which is actually the minimum Kullback-Leibler divergence estimator, the minimum Pearson chi-square estimator, the (robust) minimum Hellinger distance estimator, the minimum $L^2$ distance estimator, etc. \cite{basu2011statistical} provided a detailed discussion on a large class of MDEs along with their properties and illustrations.

There are some special classes of divergences available in the literature, which are well known for producing high robustness and efficiency properties of the resulting MDEs. One such popular class is the density power divergence (DPD) family and the associated minimum DPD estimator (MDPDE), introduced by \cite{basu1998robust},  based on a set of independent and identically distributed (IID) observations. This class of MDPDEs is seen to connect the non-robust but (asymptotically) most efficient MLE and the highly robust (but relatively less efficient) minimum $L^2$ distance estimator through a tuning parameter $\alpha$, that balances the robustness and efficiency of the MDPDEs. In recent times, these MDPDEs have become quite popular in the literature of robust inference due to their strong robustness, high efficiency, computational simplicity, and ease of interpretation. The MDPDEs have been successfully applied to many different parametric statistical models and different types of data. \cite{ghosh2013robust} extended them to the setup of independent non-homogeneous (INH) data, where the sample consists of observations that are independent but not identically distributed, that have to be modeled by a more complex parametric structure that involves a common parameter $\bm{\theta}$. This setup covers several important statistical problems, e.g., linear regression model (\cite{ghosh2013robust}), generalized linear models including logistic and Poisson regression (\cite{ghosh2016robust}), beta regression (\cite{ghosh2019robust}), non-linear regression model (\cite{jana2024robust}), etc. Most of the literature on INH data have ensured the local robustness properties of the estimators through bounded influence functions. In this paper, we shall focus on measuring the global robustness of the MDPDE for general parametric models under the setup of INH data. 

The concept of breakdown point (BP), introduced by \cite{hampel1971general}, is a very popular global measure of robustness of an estimator. In simple terms, it is defined as the minimum proportion of the sample observations that is needed to be replaced for the estimator to take any arbitrary value away from the true parameter value. \cite{basu1998robust} obtained the simultaneous asymptotic BP (ABP) of the MDPDE of the location and scale parameter of the normal distribution as $\alpha/(1+\alpha)^{3/2}$, for a given set of IID observations. Within the family of $\phi$-divergences (or disparities), \cite{park2004minimum} showed that the ABP of the corresponding MDE is at least $1/2$ at the model under a set of logical conditions, allowing only for the \say{exploding} type breakdown, where the absolute value of the estimator tends to $\infty$. Later, under the setup of INH observations \cite{ghosh2013robust} demonstrated that the ABP of the MDPDE of the location parameter of a location-scale type family is $1/2$, assuming the scale parameter to be fixed. Recently \cite{roy2023asymptotic} have derived the ABP of the minimum S-divergence estimator, a more general class of estimators containing the MDPDE, under the general parametric setup of IID observations. However, the issue of ABP of the MDPDE for general classes of parametric models beyond the IID setup remains unaddressed. Among them, the INH setup is crucial as a general framework applicable to practically all types of regression problems, statistical modeling of data with a response variable and one or more fixed covariates. Recognizing the importance of this setup, in this paper, we will derive the ABP of MDPDEs under the general INH setup and subsequently apply the general results in order to explicitly study them under some popular regression setups. We will, in fact, provide some mild conditions regarding the underlying parametric models and their behavior with the contaminating distributions, and derive some results giving the lower bound of the ABP of the minimum DPD functional (MDPDF). We will also show that the lower bound of the ABP of the MDPDE is at least as large as that of the MDPDF.\\

\noindent \textbf{Organization of the paper}: The rest of the paper is structured as follows. Section \ref{BG} provides a brief background of the DPD measure and the MDPDE under both IID and INH setups, along with mathematical definitions of the ABP of an estimator and the associated statistical functional. Our main results are contained in Section \ref{BP}, where we provide the ABP of the MDPDF and the MDPDE under suitable assumptions for a general parametric INH setup. In Section \ref{Application}, we discuss their applications to some commonly used regression setups, while Section \ref{simulation} contains extensive simulation studies, empirically illustrating the BP of the MDPDEs for these examples. Some concluding remarks are given in Section \ref{conclusion}. Finally, all the proofs and additional technical details are provided in the Appendices \ref{Proofs}-\ref{VA}. Throughout this article, we assume that the density of any probability distribution is defined with respect to the Lebesgue measure, if it is continuous, or the counting measure, if it is discrete.
% and $\mathcal{N}(\mu,\sigma^2)$ denotes the normal distribution with mean $\mu$ and variance $\sigma^2$. For any sequence of real numbers $\{a_m\}$, $a_m = O(1)$ implies that there exists $C>0$ such that $0< \liminf_{m\to\infty} |a_m| \leq \limsup_{m\to\infty} |a_m| \leq C$, and for any two sequences of real numbers $\{a_m\}, \{b_m\}$, $a_m = o(b_m)$ implies that $\lim_{m\to\infty} \frac{a_m}{b_m} = 0$

\section{Preliminary concepts} \label{BG}
\subsection{The DPD and minimum DPD estimators} \label{MDPDE-INH}
The DPD measure, first introduced by \cite{basu1998robust}, is a divergence between two densities leading to an effective robust estimation procedure. With a tuning parameter $\alpha\geq 0$, the DPD measure between two densities $g$ and $f$ (with respect to a common dominating measure $\lambda$) is given by
\begin{equation*}
    d_\alpha(g,f) =
    \begin{cases}
        \int \left\{f^{1+\alpha} - \left(1+\frac{1}{\alpha}\right)f^\alpha g + \frac{1}{\alpha}g^{1+\alpha}\right\} d\lambda, & \mbox{if }~~~ \alpha>0,\\
        \int g \ln\left(\frac{g}{f}\right) d\lambda, & \mbox{if }~~~ \alpha = 0.
    \end{cases}
\end{equation*}
Note that, by the above definition, $d_0(g,f)$ coincides with the usual Kullback-Leibler (KL) divergence, and it is obtained from the general form through its continuous limit, $\lim_{\alpha\downarrow 0}d_\alpha(g,f) = d_0(g,f)$. Also, at $\alpha=1$, it gives the squared $L^2$ distance.

Now, we consider the setup of a set of independent observations $y_1,\hdots,y_n$, with $y_i$ following a true distribution $G_i$ (unknown) with density $g_i$, for each $i=1,\hdots,n$, which we wish to model by the parametric model family $\mathcal{F}_i = \left\{f_{i,\bm{\theta}}: \bm{\theta} \in \bm{\Theta} \subseteq \mathbb{R}^d \right\}$. The model densities may be different for each $i=1,\hdots,n$, but they must depend on the same parameter $\bm{\theta}$. Then, following the approach of \cite{ghosh2013robust}, the MDPDE of $\bm{\theta}$ is defined as the minimizer of the average DPD measure between the estimated true densities and the model densities, given by $n^{-1}\sum_{i=1}^n d_\alpha(\widehat{g}_i, f_{i,\bm{\theta}})$, over the parameter space $\bm{\Theta}$, where $\widehat{g}_i$ is an empirical estimate of the true density $g_i$ based on the observation $y_i$, for each $i=1,\hdots,n$. Upon further simplification, the MDPDE $\bm{\widehat{\theta}}_\alpha$ of $\bm{\theta}$, with tuning parameter $\alpha>0$, may equivalently be obtained as the minimizer of the simpler objective function
\begin{equation}\label{dpd-obj-fn}
    H_{n,\alpha}(\bm{\theta}) = \frac{1}{n}\sum_{i=1}^n \left[\int f_{i,\bm{\theta}}^{1+\alpha}d\lambda - \left(1+\frac{1}{\alpha}\right) f_{i,\bm{\theta}}^\alpha(y_i) + \frac{1}{\alpha}\right],
\end{equation}
over $\bm{\theta\in\bm{\Theta}}$. The corresponding estimating equation is given by
\begin{equation*}
    \sum_{i=1}^n\left[\bm{u}_{i,\bm{\theta}}(y_i) f_{i,\bm{\theta}}^\alpha(y_i) - \int  \bm{u}_{i,\bm{\theta}} f_{i,\bm{\theta}}^{1 + \alpha}d\lambda\right] = 0,
\end{equation*}
where $\bm{u}_{i,\bm{\theta}}(y) = \nabla_{\bm{\theta}} \ln{f_{i,\bm{\theta}}(y)}$ is the score function of the $i$-th model density, for each $i\geq 1$ ($\nabla_{\bm{\theta}}$ denotes the partial derivative with respect to $\bm{\theta}$). Also, as $\alpha\to 0$, the objective function $H_{n,\alpha}(\bm{\theta})$ becomes (in limit) $1-n^{-1}\sum_{i=1}^n \ln{f_{i,\bm{\theta}}(y_i)}$, and hence the associated estimating equation becomes $\sum_{i=1}^n \bm{u}_{i,\bm{\theta}}(y_i) = 0$, which is the ordinary likelihood score equation. So, at $\alpha=0$, the MDPDE coincides with the classical MLE of $\bm{\theta}$. Note that the extra term $1/\alpha$ in the objective function \eqref{dpd-obj-fn} ensures the above-mentioned limiting results for $H_{n,\alpha}$ as $\alpha\rightarrow 0$.

Further, under this setup of INH data, the MDPDF of the parameter $\bm{\theta}$, at the vector of the true distributions $\bm{G} = (G_1,\hdots,G_n)$, is defined as
\begin{equation} \label{MDPDF}
    \bm{T}_\alpha(\bm{G}) = \argmin_{\bm{\theta\in\bm{\Theta}}} \frac{1}{n}\sum_{i=1}^n d_\alpha\left(g_i, f_{i,\bm{\theta}}\right) = \argmin_{\bm{\theta\in\bm{\Theta}}} H_{n,\alpha}^*(\bm{\theta}),
\end{equation}
where
\begin{equation*}
    H_{n,\alpha}^*(\bm{\theta}) = E\left(H_{n,\alpha}(\bm{\theta})\right) = \frac{1}{n}\sum_{i=1}^n \left[\int f_{i,\bm{\theta}}^{1+\alpha}d\lambda - \left(1+\frac{1}{\alpha}\right) \int f_{i,\bm{\theta}}^\alpha g_i d\lambda + \frac{1}{\alpha}\right],\ \alpha>0.
\end{equation*}
Again, this MDPDF coincides (in the limiting sense) with the maximum likelihood functional as $\alpha\rightarrow 0$. These definitions of the MDPDE and MDPDF under INH setup can be further seen to coincide with the corresponding definition under IID setup, as initially presented by \cite{basu1998robust}, when $g_i=g$ and $f_{i,\bm{\theta}}=f_{\bm{\theta}}$ for all $i\geq 1$. See \cite{ghosh2013robust} for further details on their properties and applications. In particular, the influence function (IF) of the MDPDF $\bm{T}_\alpha(\bm{G})$ was derived in Section 4 of \cite{ghosh2013robust}. The IF can be shown to be bounded whenever the term $\bm{u}_{i,\bm{\theta}}(y) f_{i,\bm{\theta}}^\alpha(y)$ is a bounded function in $y$, for each $i\geq 1$, which holds for most commonly used parametric setups for every $\alpha>0$, ensuring the local robustness of the estimator.

\subsection{Definitions of the asymptotic breakdown point} \label{BP-def}
\cite{hodges1967efficiency} first motivated the concept of a finite-sample BP of an estimator of location as the maximum proportion of contamination (erroneous values) in the sample that an estimator can tolerate without producing arbitrarily large or small values. The ABP of an estimator can be defined by taking the limit of the finite-sample BP, if it exists, as the sample size goes to $\infty$. However, this concept is restricted only to the location estimators and does not include the cases of the scale estimators, as they can break down when the estimates either \say{implode} to $0$ or \say{explode} to $\infty$ (see \cite{maronna2019robust}). Following \cite{hodges1967efficiency} and \cite{hampel1971general}, recently \cite{roy2023asymptotic} have also summarized mathematical definitions of the ABP of an estimator and the associated statistical functional under the setup of IID data. Here we extend these definitions to the setup of INH observations discussed in Section \ref{MDPDE-INH}.

Let the vector of contaminated distributions is $\bm{G}_{\epsilon,m} = \left(G_{1,\epsilon,m}, \hdots, G_{n,\epsilon,m}\right)$, where
\begin{equation*}
    G_{i,\epsilon,m} = (1-\epsilon)G_i + \epsilon K_{i,m},\ i = 1,2,\hdots,n,
\end{equation*}
with $\left\{K_{i,m}\right\}_{m\geq 1}$ being the sequences of contaminating distributions for each $i\geq 1$. Then the ABP $\epsilon^*$ of a sequence of estimators $\{\bm{T}_n\}_{n\geq 1}$ (with $n$ denoting the underlying sample size) is defined as
\begin{multline} \label{bp-def-E}
    \epsilon^* = \sup \bigg\{\epsilon\in[0,1/2]: \inf_{\bm{\theta}_\infty \in \partial\bm{\Theta}} \liminf_{n\to\infty} \liminf_{m\to\infty} P_{\bm{G}_{\epsilon,m}}\left[\lVert \bm{T}_n-\bm{\theta}_\infty\rVert > 0\right] = 1,\\
     \mbox{ for all possible contaminating sequences } \{K_{i,m}\}_{m\geq 1} \bigg\},
\end{multline}
where $P_{\bm{G}_{\epsilon,m}}$ is the probability measure under the contaminated distribution $\bm{G}_{\epsilon,m}$, $\lVert\cdot\rVert$ is the usual Euclidean distance ($L^2$-norm), and $\partial\bm{\Theta}$ is the boundary of the parameter space $\bm{\Theta} \subseteq \mathbb{R}^d$, in the extended real number system. For example, a location parameter has $\bm{\Theta} = (-\infty,\infty)$ with $\partial\bm{\Theta} = \{-\infty,\infty\}$, whereas $\bm{\Theta} = (0,\infty)$ with $\partial\bm{\Theta} = \{0, \infty\}$ for a scale parameter. Note that, the probability statement in the above definition indicates the probability that the estimator $\bm{T_n}$ is not driven to the boundary of the parameter space, under the contaminated distribution, and the ABP has been defined as the maximum proportion of contamination $\epsilon$, for which the minimum value of such probability is one.

We further define the ABP of a statistical functional $\bm{T}$, estimating a parameter $\bm{\theta}\in \bm{\Theta} \subseteq\mathbb{R}^d\ (d\geq 1)$, under the INH setup, as
\begin{multline}\label{bp-def-F}
    \epsilon^* = \sup\bigg\{\epsilon\in[0,1/2]: \inf_{\bm{\theta}_\infty \in \partial\bm{\Theta}} \liminf_{n\to \infty} \liminf_{m\to\infty} \lVert\bm{T}(\bm{G}_{\epsilon,m})-\bm{\theta}_\infty\rVert > 0,\\
    \mbox{ for all possible contaminating sequences } \{K_{i,m}\}_{m\geq 1} \bigg\}.
\end{multline}
Under the setup of IID observations, all the true distributions $G_1,G_2,\hdots$ coincide with the same underlying distribution, say $G$, and considering the same sequence of contaminating distributions, say, $\left\{K_m\right\}_{m\geq 1}$, the vector of contaminated distributions is then made of a single entity $G_{\epsilon,m} = (1-\epsilon)G + \epsilon K_m$. Then the above definitions reduce to the definitions provided in Section 2.3 of \cite{roy2023asymptotic}, under the IID setup. 

\section{Asymptotic breakdown point analysis under INH setup} \label{BP}
In this section, we will derive the ABP of the MDPDE under the setup of the INH observations for a general class of parametric model families. Let us consider the setup of the INH observations as mentioned in Section \ref{MDPDE-INH}, where $g_i$ is the true data-generating density of the observation $y_i$, which is to be modeled by the parametric model family of densities $\mathcal{F}_i$. Here we assume that the true densities $g_i$ belong to the interior of the model families $\mathcal{F}_i$, i.e., $g_i = f_{i,\bm{\theta}_0},\ i=1,\hdots,n$, for some $\bm{\theta}_0 \in \bm{\Theta} \setminus \partial\bm{\Theta}$. So, following the definition of the MDPDF from \eqref{MDPDF}, we have $\bm{\theta}_0 = \bm{T}_\alpha(\bm{G})$; ensuring its Fisher consistency. Let $k_{i,m}$ denotes the density of the contaminating distribution $K_{i,m}$ and let $g_{i,\epsilon,m} := (1-\epsilon)g_i + \epsilon k_{i,m}$ be the density of the contaminated distribution $G_{i,\epsilon,m}$, for each $i=1,\hdots,n$. Further, we denote, for any density $g$, $M_g = \int g^{1+\alpha} d\lambda$. We now present a set of assumptions, under which we shall prove our main result on the ABP of the MDPDE of $\bm{\theta}$.

\begin{enumerate}[label=(B\arabic*), ref=(B\arabic*)]
    \item\label{itm:bp1} For all $i=1,2\hdots$, the sequences of contaminating distributions $\left\{K_{i,m}\right\}_{m\geq 1}$ are such that for any compact set $S\subset\bm{\Theta}$ with $S\cap \partial\bm{\Theta}=\phi$, we have $\int \min \left\{f_{i,\bm{\theta}}, k_{i,m}\right\}d\lambda \to 0$ as $m\to \infty$, uniformly for all $\bm{\theta}\in S$ and for each $i\geq 1$.
    \item\label{itm:bp2} For any sequence of parameters $\bm{\theta}_m \to \bm{\theta}_\infty \in \partial\bm{\Theta}$ as $m\to\infty$, we have $\int\min \left\{g_i, f_{i,\bm{\theta}_m}\right\}d\lambda \to 0$ as $m\to \infty$, uniformly for each $i\geq 1$.
    \item\label{itm:bp3} The parametric model families $\mathcal{F}_i$ and the families of contaminating densities $\{k_{i,m}\}_{m\geq 1},\ i=1,2,\hdots$, are uniformly $L^{1+\alpha}$ integrable, so that
    \begin{equation*}
    \sup_{i\in \mathbb{N}} \limsup_{m\to\infty} \int k_{i,m}^{1+\alpha}d\lambda < \infty, \mbox{ and }~~ \sup_{i\in \mathbb{N}} \sup_{\bm{\theta\in\Theta}} \int f_{i,\bm{\theta}}^{1+\alpha}d\lambda < \infty.
    \end{equation*}
    \item\label{itm:bp4} There exists $\epsilon_\alpha \in \left(0,\frac{1}{2}\right]$ such that for all $\epsilon \in \left[0, \epsilon_\alpha\right)$ and for sufficiently large $n$,
    \begin{equation}\label{bp3ineq}
        \liminf_{m\to\infty} \frac{1}{n}\sum_{i=1}^n d_\alpha(\epsilon k_{i,m}, f_{i,\bm{\theta}_m}) > \frac{\epsilon^{1+\alpha}}{\alpha}\limsup_{m\to \infty} \frac{1}{n}\sum_{i=1}^n M_{k_{i,m}} + q_\alpha(1-\epsilon) \frac{1}{n}\sum_{i=1}^n M_{g_i},
    \end{equation}
    for any sequence of parameters $\bm{\theta}_m \to \bm{\theta}_\infty \in \partial\bm{\Theta}$ as $m\to\infty$, where $q_\alpha(\epsilon) = 1 - \frac{1+\alpha}{\alpha}\epsilon$, for any $\epsilon \in (0,1]$.
\end{enumerate}

Note that Assumption \ref{itm:bp1} ensures that the sequences of contaminating distributions $\left\{K_{i,m}\right\}_{m\geq 1}$ are asymptotically singular to the true distributions and to the specific model distributions, for each $i=1,2,\hdots$, within a compact subset of the parametric family. Assumption \ref{itm:bp2} ensures that the true distributions are asymptotically singular to the model distributions with parameter tending to the boundary of the parameter space, for each $i\geq 1$. Assumption \ref{itm:bp4} is the most important assumption, providing the required lower bound of the breakdown point of interest, and \ref{itm:bp3} is just a technical requirement to ensure boundedness of certain quantities. The following theorem presents formal results on the ABP of the MDPDF of $\bm{\theta}$, with $\alpha>0$, under the INH setup. Its proof is given in Appendix \ref{pf-bp-th}. Note that, at $\alpha=0$, the MDPDF coincides with the maximum likelihood functional, having the ABP zero under most parametric setups, and hence we consider only the MDPDFs with $\alpha>0$ in the following theorem.

\begin{theorem}\label{bp-th}
    Under Assumptions \ref{itm:bp1}-\ref{itm:bp4}, if the true densities $g_i$ belong to the interior of the model families $\mathcal{F}_i$ i.e., $g_i = f_{i,\bm{\theta}_0},\ i=1,\hdots,n$, for some $\bm{\theta}_0 \in \bm{\Theta} \setminus \partial\bm{\Theta}$, then, for any fixed $\alpha>0$, the asymptotic breakdown point of the MDPDF $\bm{T}_\alpha$ is at least $\epsilon_\alpha$, where $\epsilon_\alpha$ is as defined in Assumption \ref{itm:bp4}.
\end{theorem}

The above theorem does not depend on any special property of the boundary of the parameter space $\partial\bm{\Theta}$, except the singularity Assumptions \ref{itm:bp1} and \ref{itm:bp2}. So, it remains applicable even if these assumptions hold for some specific points $\bm{\theta}_\infty \in \partial\bm{\Theta}$ or a subset of $\partial\bm{\Theta}$. It allows us to study specific types of breakdown separately by verifying these Assumptions \ref{itm:bp1}-\ref{itm:bp2} only for the (sub)set of boundary points of interest. For example, when the parameter space of a scale parameter $\sigma$ is $\bm{\Theta}_{\sigma}=(0,\infty)$, then its estimator can possess \say{imploding} or \say{exploding}-type of breakdown at $\bm{\theta}_\infty=0$ or $\bm{\theta}_\infty=\infty$, respectively. However, the \say{exploding} type breakdown can only be considered by applying the above theorem to the subset $\{+\infty\}$ of $\partial \bm{\Theta}_\sigma$.

However, in practice, verifying Assumption \ref{itm:bp4} and deriving the lower bound $\epsilon_\alpha$ of the ABP (as in Theorem \ref{bp-th}) under different parametric setups may be difficult. So, we present three alternative assumptions, which are individually sufficient to ensure Assumption \ref{itm:bp4} and comparatively easier to verify and use in specific cases.

\begin{enumerate}[label=(B5), ref=(B5)]
    \item\label{itm:bp5} $M_{k_{i,m}}\leq M_{f_{i,\bm{\theta}_m}}$, for all sufficiently large $m$ and for each $i\geq 1$.
\end{enumerate}

Based on an application of the Corollary 3.1 of \cite{roy2023asymptotic}, Assumption \ref{itm:bp5} implies Assumption \ref{itm:bp4} with $\epsilon_\alpha = \alpha/(1+\alpha)$. So, in Theorem \ref{bp-th}, one can replace Assumption \ref{itm:bp4} with Assumption \ref{itm:bp5} to deduce the following corollary. Its formal proof is given in Appendix \ref{bp-th-2-pf}.
\begin{corollary} \label{bp-th-2}
    Under Assumptions \ref{itm:bp1}-\ref{itm:bp3} and \ref{itm:bp5}, if the true densities $g_i$ belong to the interior of the model families $\mathcal{F}_i$, i.e., $g_i = f_{i,\bm{\theta}_0},\ i=1,\hdots,n$, for some $\bm{\theta}_0 \in \bm{\Theta} \setminus \partial\bm{\Theta}$, then the asymptotic breakdown point of the MDPDF $\bm{T}_\alpha$ is at least $\epsilon_\alpha^* = \alpha/(1+\alpha)$ for any given $\alpha\in[0,1]$.
\end{corollary}

\begin{example} \label{ex-bp-th-2}
    Given a set of $n$ INH observations, consider the problem of estimating the location parameter of the location families of densities $\mathcal{F}_i=\left\{f_{i,\bm{\theta}}: \bm{\theta}\in\bm{\Theta} \subseteq\mathbb{R}^d \right\}$, where $f_{i,\bm{\theta}}(x) = f\left(x-l_i(\bm{\theta})\right) \mbox{ for all } x\in\mathbb{R}$, $l_i:\mathbb{R}^d \to \mathbb{R}$ are some one to one functions for each $i=1,\hdots,n$, and $f$ is a known probability density (e.g., standard normal density). Also, suppose that, for a fixed $\bm{\theta}\in\bm{\Theta}$, $|l_i(\bm{\theta})|$ is bounded in $i$ and breakdown occurs for a sequence of parameters $\bm{\theta}_m$ if $|l_i(\bm{\theta}_m)|\to\infty, \mbox{ as } m\to\infty$, for each $i\geq 1$. We assume that $g_i = f_{i,\bm{\theta}_0}$, for $i=1,\hdots,n$, and some $\bm{\theta}_0 \in \bm{\Theta} \setminus \partial\bm{\Theta}$, and that the location families ($\mathcal{F}_i$) are such that Assumption \ref{itm:bp2} is satisfied and $\sup_{i\geq 1} \sup_{\bm{\theta\in\Theta}} \int f_{i,\bm{\theta}}^{1+\alpha} d\lambda < \infty$. Then, considering the sequences of contaminating densities $\left\{k_{i,m} \right\}_{m\geq 1}$ from $\mathcal{F}_i$, with $k_{i,m} \equiv f_{i,\bm{\theta}^*_m}$, where $|l_i(\bm{\theta}^*_m)|\to\infty \mbox{ as } m\to\infty$, for each $i=1,\hdots,n$, Assumptions \ref{itm:bp1}, \ref{itm:bp3} and \ref{itm:bp5} are satisfied. Now, applying Corollary \ref{bp-th-2}, we get the ABP of the MDPDF to be at least $\alpha/(1+\alpha)$, for any $\alpha\in[0,1]$. However, a stricter lower bound can be obtained by directly applying Theorem \ref{bp-th}. In the Inequality \eqref{bp3ineq} of Assumption \ref{itm:bp4}, an application of H\"{o}lder's inequality to the term $\int f_{i,\bm{\theta}_m}^\alpha k_{i,m} d\lambda$ within $\frac{1}{n}\sum_{i=1}^n d_\alpha(\epsilon k_{i,m}, f_{i,\bm{\theta}_m})$ yields that a sufficient condition for Inequality \eqref{bp3ineq} to hold is $\epsilon<1/2$. So, the lower bound of the ABP is $1/2$ for the MDPDF of the location parameter when both parametric models and contaminating distributions come from the same location family for each $i\geq 1$. This result was first obtained by \cite{ghosh2013robust} with applications to a linear regression model with IID and normally distributed errors having known variance. \qed
\end{example}

We have seen, in the above example, that Assumption \ref{itm:bp5} is much stricter than Assumption \ref{itm:bp4}, albeit being easier to verify. We now present another alternative sufficient condition for Assumption \ref{itm:bp4} that can also be easily verified only with the knowledge of the sequences of contaminating densities $\{k_{i,m}\}_{m\geq 1}$. One can check the validity of this assumption to ensure \ref{itm:bp4}, even if \ref{itm:bp5} is not satisfied or is difficult to verify.

\begin{enumerate}[label=(B6), ref=(B6)]
    \item\label{itm:bp6} There exists a constant $C\in [0,\infty)$ such that $\limsup_{m\to \infty}M_{k_{i,m}} \leq C$, for each $i\geq 1$.
\end{enumerate}

Under Assumption \ref{itm:bp6}, the following corollary ensures the validity of Assumption \ref{itm:bp4} and provides an easily obtainable lower bound of the breakdown point of interest. Its proof is given in Appendix \ref{bp-th-3-pf}.
\begin{corollary}\label{bp-th-3}
    Under Assumptions \ref{itm:bp1}-\ref{itm:bp3} and \ref{itm:bp6}, if the true densities $g_i$ belong to the interior of the model families $\mathcal{F}_i$, i.e., $g_i = f_{i,\bm{\theta}_0},\ i=1,\hdots,n$, for some $\bm{\theta}_0 \in \bm{\Theta} \setminus \partial\bm{\Theta}$, then, for any $\alpha>0$, the asymptotic breakdown point of the MDPDF $\bm{T}_\alpha$ is at least $\min\{\epsilon_\alpha^\prime, \frac{1}{2}\}$, where $\epsilon_\alpha^\prime$ is the unique solution, in the interval $(0,1)$, of the equation
    \begin{equation}\label{cor-eq}
        \frac{C}{\alpha}x^{1+\alpha} + q_\alpha(1-x) L_0 = 0, \mbox{ where } L_0 = \lim_{n\to\infty}\frac{1}{n}\sum_{i=1}^n M_{g_i},
    \end{equation}
    and the function $q_\alpha(\cdot)$ is as defined in Assumption \ref{itm:bp4}.
\end{corollary}
\begin{remark} \label{remark-bp-th-3}
    For a certain problem, Assumption \ref{itm:bp6} holds with $C=0$. Then it is easy to observe that Eq.~\eqref{cor-eq} of Corollary \ref{bp-th-3} has the solution $1/(1+\alpha)$ and hence a lower bound of the ABP is $\min \left\{\dfrac{1}{1+\alpha}, \dfrac{1}{2}\right\}=\dfrac{1}{2}$, for $0<\alpha\leq1$.
\end{remark}
If the contaminating distributions are discrete with the probability mass functions (PMF) $k_{i,m}$, then Assumption \ref{itm:bp6} always holds with $C=1$. However, it may not be a sharp upper bound to the inequality in \ref{itm:bp6}. We should choose a sharper upper bound in \ref{itm:bp6} to get a sharper lower bound of the ABP of the MDPDF through Corollary \ref{bp-th-3}.
\begin{example}
    Extending Example \ref{ex-bp-th-2}, let us consider the location-scale type parametric family $\mathcal{F}_i$, given by the collection of densities $f_{i,\bm{\theta}} = \frac{1}{\sigma} f\left(\frac{x-l_i(\bm{\mu})}{\sigma}\right)$, where $\bm{\theta}=(\bm{\mu}^T,\sigma)^T\in \bm{\Theta}\subseteq \mathbb{R}^p\times\mathbb{R}^+$. Here we consider the breakdown of a sequence of mean parameters $\bm{\mu}_m$ as $|l_i(\bm{\mu}_m)|\to\infty \mbox{ with } m\to\infty$, for each $i$, and an \say{exploding}-type breakdown for any sequence of scale parameter as $\sigma_m\to\infty$ with $m\to\infty$. Then \ref{itm:bp1}-\ref{itm:bp3} and \ref{itm:bp6} are satisfied with $C=0$, and hence by Remark \ref{remark-bp-th-3}, the ABP of the MDPDF is at least $1/2$ for any $\alpha\in (0,1]$. \qed
\end{example}

Although useful in commonly used cases, Assumptions \ref{itm:bp5} or \ref{itm:bp6} may not provide a very precise lower bound to the ABPs in some cases since they cover all uniformly integrable contaminating densities. Also, $M_{g_i}$s may be unknown, which is essential to compute the ABP through Corollary \ref{bp-th-3}, based on \ref{itm:bp6}, unless $C=0$. So, we now present another sufficient assumption for \ref{itm:bp4} utilizing the orthogonality between the sequences of model densities $\left\{f_{i,\bm{\theta}_m}\right\}_{m\geq 1}$ and contaminating densities $\left\{k_{i,m}\right\}_{m\geq 1}$; it basically says that $\sum_{i=1}^n \int f_{i,\bm{\theta}_m}^\alpha k_{i,m} d\lambda$ converges to $0$ faster than $\sum_{i=1}^n M_{f_{i,\bm{\theta}_m}}$ in a suitable asymptotic sense.
\begin{enumerate}[label=(B7), ref=(B7)]
    \item\label{itm:bp7} The sequences of contaminating densities $\left\{k_{i,m}\right\}_{m\geq 1}$, for each $i=1,\hdots,n$, satisfy
    \begin{equation*}
        L := \liminf_{\substack{m\to\infty,\\ \bm{\theta}_m \to \bm{\theta}_\infty \in \bm{\Theta}}} L_m > 0, \mbox{ where } L_m = \lim_{n\to\infty} \frac{\sum_{i=1}^n M_{f_{i,\bm{\theta}_m}}}{\sum_{i=1}^n \int f_{i,\bm{\theta}_m}^\alpha k_{i,m}d\lambda}.
    \end{equation*}
\end{enumerate}
\begin{corollary} \label{bp-th-4}
    Under Assumptions \ref{itm:bp1}-\ref{itm:bp3} and \ref{itm:bp7}, if the true densities $g_i$ belong to the interior of the model families $\mathcal{F}_i$, i.e., $g_i = f_{i,\bm{\theta}_0},\ i=1,\hdots,n$, for some $\bm{\theta}_0 \in \bm{\Theta} \setminus \partial\bm{\Theta}$, then the asymptotic breakdown point of the MDPDF $\bm{T}_\alpha$ is at least $\min\left\{ \frac{\alpha L}{1+\alpha}, \frac{1}{2}\right\}$ for any $\alpha\in [0,1]$, where $L$ is as given in Assumption \ref{itm:bp7}.
\end{corollary}
\noindent Its proof is given in Appendix \ref{bp-th-4-pf}.\\

Till now, we have discussed the ABP of the MDPDF $\bm{T}_\alpha$, where the true distributions in $\bm{G}=(G_1,\hdots,G_n)$ are assumed to be unknown. But, in reality, $\bm{G}$ needs to be replaced by its empirical estimate $\bm{\widehat{G}} = (\widehat{G}_1,\hdots,\widehat{G}_n)$ (where $\widehat{G}_i$ be the empirical distribution function associated with observation $y_i$) to obtain the MDPDE $\bm{\widehat{\theta}}_\alpha = \bm{T}_\alpha(\bm{\widehat{G}}) = \argmin_{\bm{\theta}\in\bm{\Theta}} H_{\alpha,n}(\bm{\theta})$, which has been discussed in Section \ref{MDPDE-INH}. Our next theorem, under certain assumptions, ensures that the ABP of the MDPDE $\bm{T}_\alpha(\bm{\widehat{G}})$ is at least as large as that of the associated functional $\bm{T}_\alpha(\bm{G})$ for any $\alpha>0$. Its proof is given in Appendix \ref{pf-bp-th-5}.

% \begin{theorem} \label{bp-th-5}
%     Suppose $\alpha>0$ and the following assumptions hold along with Assumption \ref{itm:bp3}.
%     \begin{enumerate}[label=(\roman*)]
%         \item With probability tending to 1 and for all $\alpha \in [0,1]$
%               \[\frac{1}{n}\sum_{i=1}^n \left(\int\left|\widehat{g}_i - g_i\right|^{\alpha+1} \right)^{\frac{1}{\alpha+1}} \to 0, \mbox{ as } n\to\infty,\]
%         \item $\bm{T}_\alpha(\bm{G})$ is a well-separated minimizer of the average DPD measure $\frac{1}{n}\sum_{i=1}^n d_\alpha\left(g_i, f_{i,\bm{\theta}}\right)$, i.e., for any $\xi > 0$, independent of $g_i$'s,
%         \[\frac{1}{n}\sum_{i=1}^n d_\alpha\left(g_i,f_{i,\bm{T}_\alpha(\bm{G})}\right) < \inf_{|\bm{\theta}-\bm{T}_\alpha(\bm{G})|>\xi} \frac{1}{n}\sum_{i=1}^n d_\alpha\left(g_i, f_{i,\bm{\theta}}\right), \mbox{ for every n}.\]
%     \end{enumerate}
%     Now, if $\epsilon_1^*,\epsilon_2^* \in [0,1/2)$ are the ABPs of the MDPDF $\bm{T}_\alpha(\bm{G})$ as given in \eqref{bp-def-F} and the MDPDE $\bm{T}_\alpha\left(\bm{\widehat{G}}\right)$, as given in \eqref{bp-def-E}, respectively, then we must have $\epsilon_2^* \geq \epsilon_1^*$.
% \end{theorem}
\begin{theorem} \label{bp-th-5}
    Suppose the Assumptions (A1)-(A6) of \cite{ghosh2013robust} hold. Further suppose that the sequences contaminating distributions $\{K_{i,m}\}_{m\geq 1}$, associated with the contaminated distributions $G_{i,\epsilon,m}$ are same for each $i\geq 1$, and is denoted by $\{K_m\}_{m\geq 1}$, with $K_m\in \mathcal{K}$ satisfy Assumption~\ref{itm:bp3} with $(\mathcal{K}, \Vert \cdot\Vert_{TV})$ being compact, and $\bm{u}_{i,\bm{\theta}}(y) f_{i,\bm{\theta}}^\alpha(y)$ is a bounded function in $y$, for each $i\geq 1$. Now, if $\epsilon_1^*,\epsilon_2^* \in [0,1/2)$ are the ABPs of the MDPDF $\bm{T}_\alpha(\bm{G})$ as defined in \eqref{bp-def-F} and the MDPDE $\bm{T}_\alpha(\bm{\widehat{G}})$ as defined in \eqref{bp-def-E}, respectively, with same tuning parameter $\alpha>0$, then we must have $\epsilon_2^* \geq \epsilon_1^*$.
\end{theorem}

In the statement of the above theorem, $\Vert \cdot\Vert_{TV}$ denotes the total variation norm. Note that the above theorem requires the parametric families to satisfy the boundedness condition of the functions $\bm{u}_{i,\bm{\theta}}(y) f_{i,\bm{\theta}}^\alpha(y)$, for all $i$. This condition also ensures the boundedness of the influence function of the MDPDE, as discussed in Section 4 of \cite{ghosh2013robust}. In other words, this boundedness condition is important to achieve both local and global robustness of the MDPDEs, which is intuitively expected given the form of its estimating equation.

\section{Applications to fixed-design regression setups} \label{Application}
We now study the ABP of the MDPDF under some useful fixed-design regression models based on the general theory developed in the preceding section. Let us start by stating two conditions regarding any double-sequence of real numbers $\{a_{i,m}\}_{i\geq1, m\geq 1}$, which will be used throughout the rest of the paper.
\begin{enumerate}[label=(C\arabic*), ref=(C\arabic*)]
\item \label{itm:c1} There exists a lower bound $l_m$ and an upper bound $u_m$ with $l_m\leq a_{i,m}\leq u_m$, for all $i\geq 1$ and $m\geq1$.
\item \label{itm:c2} There exists a lower bound $l_m$ and an upper bound $u_m$ with $l_m\leq a_{i,m}\leq u_m$, for all $i\geq 1,~m\geq1$, and both $l_m$ and $u_m$ tend to $\infty$ as $m\to\infty$ and $l_m/u_m = O(1)$
\end{enumerate}

\subsection{Linear and non-linear regression with normally distributed errors} \label{RegModels}
Consider a typical non-linear regression (NLR) model of the form
\begin{equation} \label{nlr}
    y_i = \mu(\bm{x}_i,\bm{\beta}) + \epsilon_i,\ i=1,\hdots,n,
\end{equation}
where $y_i$ and $\bm{x}_i$ are, respectively, the values of the response variable $Y$ and covariates $X_1,\hdots,X_p$, for the $i$-th observation, $\bm{\beta}$ is a $k$-dimensional ($k$ could be different from $p$) regression parameter vector, and the random errors $\epsilon_i$s are assumed to be the independently distributed as normals with common mean $0$ and common variance $\sigma^2$. Here, the mean function $\mu(\bm{x_i},\bm{\beta})$ is generally a known function in the (unknown) parameter $\bm{\beta}$. However, as a special case, \eqref{nlr} leads to the standard linear regression model at the choice $\mu(\bm{x}_i, \bm{\beta})=\bm{x}_i^T\bm{\beta}$. The MDPDE of the parameters of a general NLR model, as in \eqref{nlr}, has been studied in \cite{jana2024robust}, along with its statistical properties. Here, the assumed parametric model family of $y_i$ is 
\[\mathcal{F}_i = \left\{\mathcal{N}\left(\mu(\bm{x_i},\bm{\beta}), \sigma^2 \right): \bm{\beta}\in\bm{\Theta}_{\bm{\beta}} \subseteq\mathbb{R}^{k}, \sigma > 0 \right\}.\]
Let us denote $\bm{\theta} = (\bm{\beta}^T,\sigma)^T,\ \bm{\Theta} = \bm{\Theta}_{\bm{\beta}} \times(0,\infty)$ and $f_{i,\bm{\theta}}$ be the density of $\mathcal{N}\left(\mu(\bm{x_i},\bm{\beta}), \sigma^2 \right)$. We also assume that the true density of $y_i$ belongs to the corresponding model family $\mathcal{F}_i$, for each $i=1,\hdots,n$, with true parameter values being $\bm{\theta}_0 = (\bm{\beta}_0^T, \sigma_0)^T \in \bm{\Theta}\setminus \partial\bm{\Theta}$.

Now, we consider a sequence of parameter vectors $\bm{\theta}_m = (\bm{\beta}_m^T, \sigma_m)^T$, for $m\geq 1$. At any $\alpha>0$, we have $M_{f_{i,\bm{\theta}_m}} = \sigma_m^{-\alpha} M_\phi$ (with $\phi$ being the standard normal density), which does not remain bounded if $\sigma_m \to 0$ as $m \to \infty$. So, to avoid this problem, we only consider the \say{exploding}-type breakdown for $\sigma$. In other words, we consider \sr{the simultaneous breakdown of $\bm{\theta}$ such that $\bm{\theta}_m := (\bm{\beta}_m^T, \sigma_m)^T$ satisfy $\sigma_m \to \infty$ and $\vert \mu(\bm{x}_i, \bm{\beta}_m)\vert \to \infty$ as $m \to \infty$ for all $i$.} But, for a fixed $\bm{\beta}\in\bm{\Theta}_{\bm{\beta}}$, we assume that $\mu(\bm{x_i},\bm{\beta})$ remains uniformly bounded in $i$, i.e., there exists two constants, $l_{\bm{\beta}},~u_{\bm{\beta}}\in \mathbb{R}$ (independent of $i$) such that $l_{\bm{\beta}}\leq \mu(\bm{x_i},\bm{\beta}) \leq u_{\bm{\beta}}$, for all $i\geq 1$. 

\sr{Let us consider the specific case when the contaminating distributions $K_{i,m}$ are given by $\mathcal{N} \left(\eta_{i,m}, \delta_{i,m}^2\right)$ with $|\eta_{i,m}|\to\infty$ and $\delta_{i,m} \to \infty$, as $m\to\infty$ for each $i\geq 1$ and $\eta_{i,m}$ satisfies Condition \ref{itm:c1}. Then it can be verified that Assumptions~\ref{itm:bp1}-\ref{itm:bp3} hold. Additionally, the contaminating densities $k_{i,m}$ satisfy~\ref{itm:bp5}, i.e, $M_{k_{i,m}} \to 0$, as $m \to \infty$. Therefore, we can apply Corollary \ref{bp-th-2} to have the ABP of the MDPDF to be at least $\epsilon_\alpha^* = \alpha/(1+\alpha)$ for any $\alpha\in[0,1]$.}

However, under the present situation, Assumption \ref{itm:bp6} is also satisfied with $C=0$. Hence, Corollary \ref{bp-th-3} and Remark \ref{remark-bp-th-3} can be used to derive that the ABP is at least $1/2$, for any $\alpha\in (0,1]$, \sr{which ensures the highest possible breakdown for the MDPDF}.

Now, we shall empirically illustrate the simultaneous BP of the MDPDF of the NLR model parameter $\bm{\theta} = (\bm{\beta}^T, \sigma)^T$ for two cases: a linear regression model and a Michaelis-Menten model (see, e.g., \cite{jana2024robust, michaelis1913kinetik}). For this purpose, let us first derive the expression of MDPDF of the parameter $\bm{\theta}$ \sr{for} a general NLR model \eqref{nlr}, \sr{when the observations in the sample are distributed according to $\bm{G}_{\epsilon,m}$}, and the contaminating distributions $K_{i,m}$ are taken as $\mathcal{N}(\mu(\bm{x_i},\bm{\beta}_m), \sigma_m^2)$, for each $i=1,\hdots,n$. Denoting the density of $\mathcal{N}(\mu,\sigma^2)$ by $\phi_{\mu,\sigma}(x) = \frac{1}{\sigma} \phi\left(\frac{x-\mu}{\sigma}\right)$, we can deduce that (see, e.g., supplementary material of \cite{roy2023asymptotic})
\begin{equation} \label{pdf-int-roy}
     \int \phi_{\mu_1,\sigma_1}(x) \phi_{\mu_2,\sigma_2}^\alpha(x)dx = (2\pi)^{-\alpha/2} \frac{\sigma_2^{1-\alpha}}{\sqrt{\alpha\sigma_1^2 + \sigma_2^2}} \exp\left\{-\frac{\alpha(\mu_1-\mu_2)^2}{2\left(\alpha\sigma_1^2 + \sigma_2^2 \right)}\right\}.
\end{equation}
Further, for simplicity of notation, let us denote $\mu_i = \mu_i(\bm{\beta}) = \mu(\bm{x_i},\bm{\beta})$, $\mu_{0i} = \mu(\bm{x_i},\bm{\beta_0})$ and $\mu_{mi} = \mu(\bm{x_i},\bm{\beta_m})$, for each $i=1,\hdots,n$. Then, for the $i$-th observation, the model density is $f_{i,\bm{\theta}} \equiv \phi_{\mu_i,\sigma}$ and the density of the contaminated distribution is $g_{i,\epsilon,m} = (1-\epsilon) \phi_{\mu_{0i},\sigma_0} + \epsilon \phi_{\mu_{mi},\sigma_m}$. But, from the definition \eqref{MDPDF}, with $f_{i,\bm{\theta}} = \phi_{\mu_i,\sigma}$ and $g_i=g_{i,\epsilon,m}$ for $i\geq 1$, an application of \eqref{pdf-int-roy} suggests that the objective function $H_{n,\alpha}^*(\bm{\theta})$ of the MDPDFs, at the contaminated distributions $g_{i,\epsilon,m}$, $i\geq 1$, for $\alpha>0$, can be simplified to

\begin{multline} \label{MDPDFnlr}
    H_{n,\alpha}^*(\bm{\theta}) = \frac{1}{\sigma^\alpha \sqrt{1+\alpha}} - \left(1 + \frac{1}{\alpha}\right)\sigma^{1-\alpha}
    \frac{1}{n}\sum_{i=1}^n \left[ \frac{1-\epsilon}{\sqrt{\alpha\sigma_0^2 + \sigma^2}}\exp\left\{-\frac{\alpha(\mu_i-\mu_{0i})^2}{2\left(\alpha\sigma_0^2 + \sigma^2\right)}\right\}\right.\\ 
    + \left.\frac{\epsilon}{\sqrt{\alpha\sigma_m^2 + \sigma^2}}\exp\left\{-\frac{\alpha(\mu_i-\mu_{mi})^2}{2\left(\alpha\sigma_m^2 + \sigma^2\right)}\right\}\right] + \frac{1}{\alpha}.
\end{multline}

\begin{example}
Following the NLR model \eqref{nlr}, we first consider a simple linear regression (SLR) model with the mean function given by
\begin{equation} \label{lr}
    \mu(\bm{x_i},\bm{\beta}) = \beta_0 + \beta_1 x_i,\ \ \ i=1,\hdots,n,
\end{equation}
where $x_i$ is the value of a scalar covariate for the $i$-th observation, for each $i=1,\hdots,n$ and $\beta_0, \beta_1$ are the regression parameters. The covariates were taken as $20$ randomly generated observations from $\mathcal{N}(50,20^2)$, the true parameter value as $\bm{\theta}_0 = \left(\bm{\beta}_0^T, \sigma_0\right)^T = (35,1,1.2)$. For the purpose of numerical illustration, fixed contaminating distributions were chosen from the model family with parameters $(\bm{\beta}^{*T}, \sigma^*)^T = (50,2,0.5)$, instead of taking contaminating sequences. We computed the values of the MDPDFs of the model parameters $\beta_0, \beta_1, \sigma$ at various contamination proportions $\epsilon$, by minimizing $H_{n,\alpha}^*(\bm{\theta})$ from Equation \eqref{MDPDFnlr} for this specific model given in \eqref{lr}. The resulting values of the MDPDFs are plotted in Figure \ref{fig:bp-lr}, for various values of $\alpha\geq 0$.
\begin{figure}[H]
     \centering
     \begin{subfigure}[b]{0.45\textwidth}
         \centering
         \includegraphics[width=\textwidth]{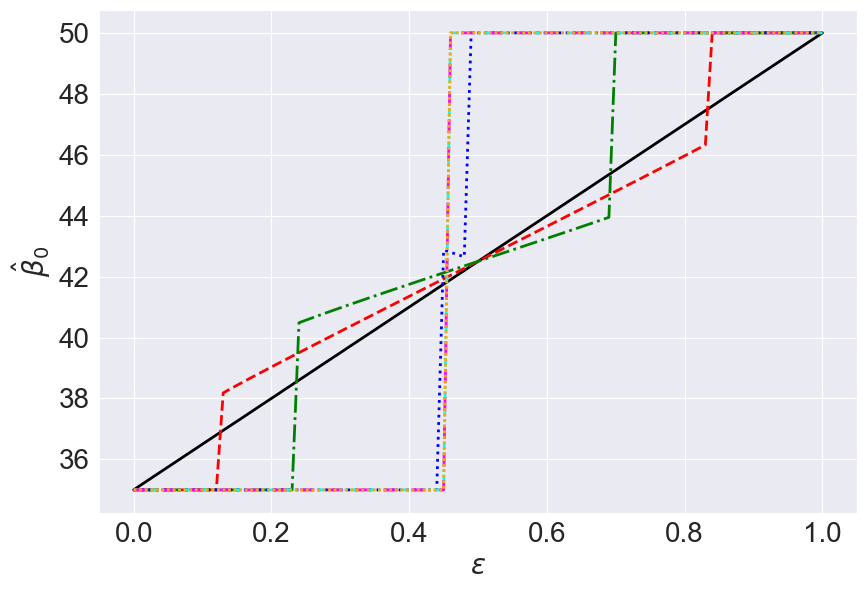}
         \caption{MDPDF of $\beta_0$}
         \label{fig:bp-lr-b0}
     \end{subfigure}
     \hfill
     \begin{subfigure}[b]{0.45\textwidth}
         \centering
         \includegraphics[width=\textwidth]{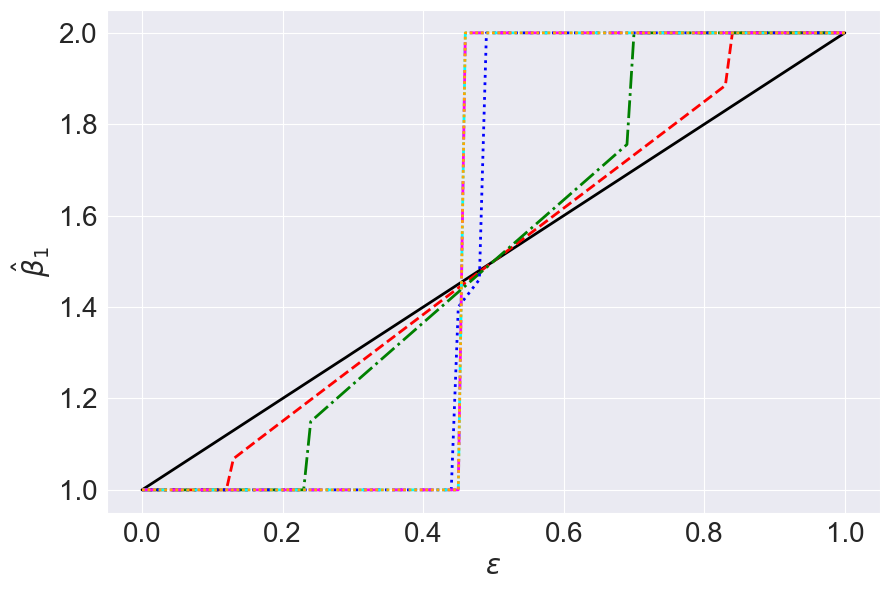}
         \caption{MDPDF of $\beta_1$}
         \label{fig:bp-lr-b1}
     \end{subfigure}
     \hfill
     \begin{subfigure}[b]{0.45\textwidth}
         \centering
         \includegraphics[width=\textwidth]{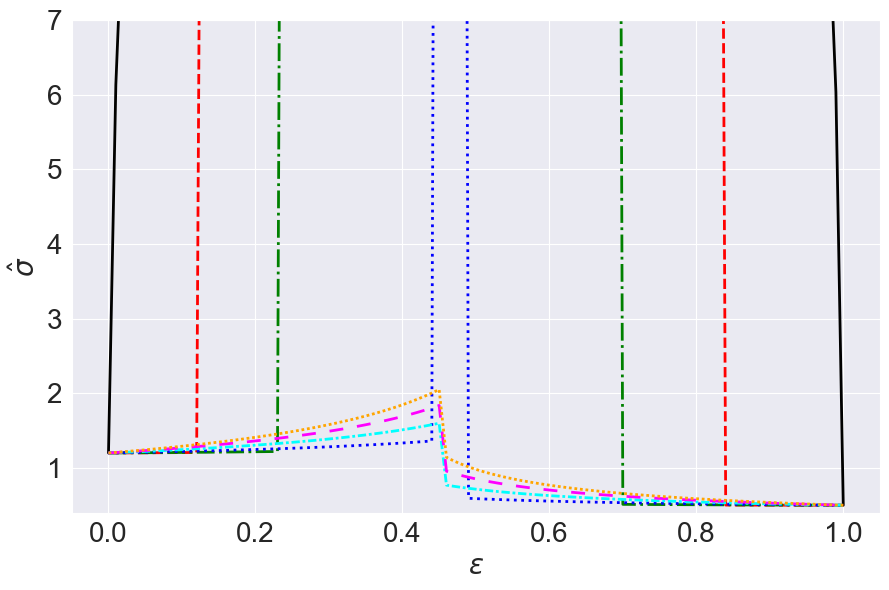}
         \caption{MDPDF of $\sigma$}
         \label{fig:bp-lr-sig}
     \end{subfigure}
     \hfill
     \begin{subfigure}[b]{0.45\textwidth}
         \centering
         \includegraphics[width=0.4\textwidth]{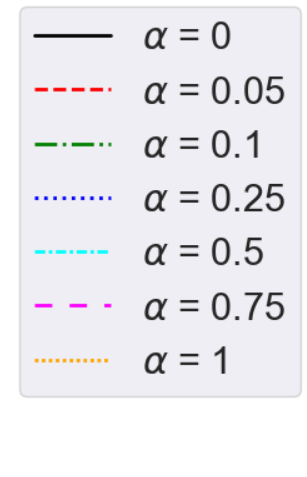}
         \caption*{}
     \end{subfigure}
     \caption{The MDPDFs obtained under the SLR model against various contamination proportions (with a fixed contamination from the model family)}
     \label{fig:bp-lr}
\end{figure}

We can clearly see from Figure \ref{fig:bp-lr} that the MDPDF with $\alpha=0$ starts to break down as soon as contamination is introduced in the underlying distributions, and the BP increases with an increase in $\alpha$. MDPDFs of both regression parameters ($\beta_0,\beta_1$) and the error standard deviation ($\sigma$) are seen to have the BPs near 0.13 and 0.24, for $\alpha=0.05$ and $0.1$, respectively. However, for all $\alpha\geq 0.25$, the BPs are observed near 0.45. For lower values of $\alpha$, the MDPDFs of $\sigma$ show an interesting pattern--they break down after a certain contamination proportion, then start increasing, reach a maxima, then again decrease, and finally coincide with the $\sigma$-value of the contaminating distribution, and thereby creating a half-circular pattern. However, we have truncated the y-axis of Fig. \ref{fig:bp-lr-sig} to maintain proper visibility of the breakdown behavior of the curves with higher $\alpha$s. \qed
\end{example}

\begin{example}
As a genuine NLR model, we consider the popular Michaelis-Menten (MM) model, which is widely used in enzyme kinetics to model the velocity ($y_i$) of an enzyme kinetic reaction with the substrate concentration ($x_i$). Its mean function is given by
\begin{equation} \label{mmnlr}
    \mu(x_i,\bm{\beta}) = \frac{\beta_1 x_i}{\beta_2 + x_i},\ i=1,\hdots,n,
\end{equation}
where the parameters $\beta_1, \beta_2$ denote the `initial reaction rate' and `Michaelis-constant', respectively; see \cite{jana2024robust} for the study of the MDPDE of these parameters under the MM model. Here, we obtained the MDPDFs of the parameters $\beta_1, \beta_2, \sigma$, considering the true parameter values to be $\left(\bm{\beta}_0^T, \sigma_0\right) = (5,2,0.5)$ and the parameters of the fixed contaminating distributions (from the model family) to be $\left(\bm{\beta}^{*T}, \sigma^*\right) = (20,3,0.1)$, with the covariate values being $n=80$ equidistant points from $0.1$ to $80$. The resulting MDPDF values are shown in Figure \ref{fig:bp-mm} for various $\alpha\geq 0$.

\begin{figure}[H]
     \centering
     \begin{subfigure}[b]{0.45\textwidth}
         \centering
         \includegraphics[width=\textwidth]{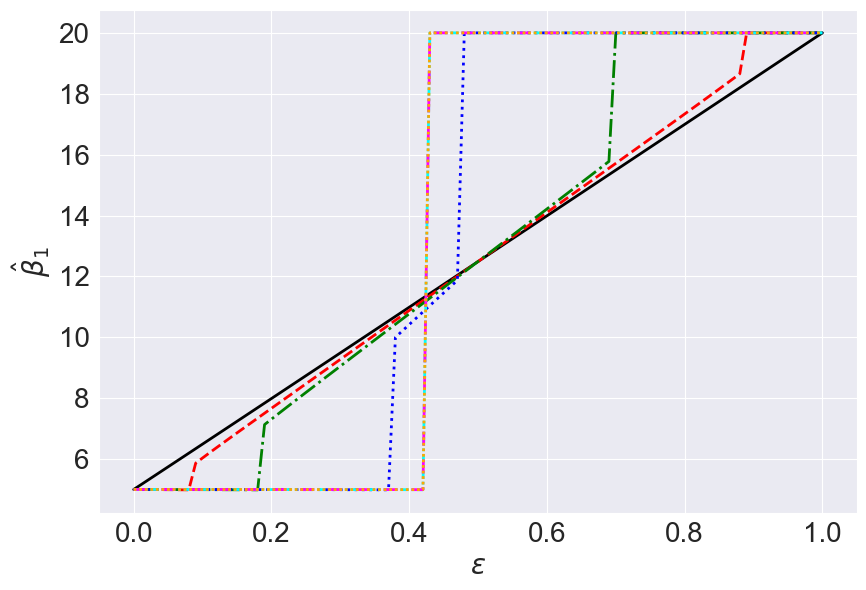}
         \caption{MDPDF of $\beta_1$}
         \label{fig:bp-mm-b1}
     \end{subfigure}
     \hfill
     \begin{subfigure}[b]{0.45\textwidth}
         \centering
         \includegraphics[width=\textwidth]{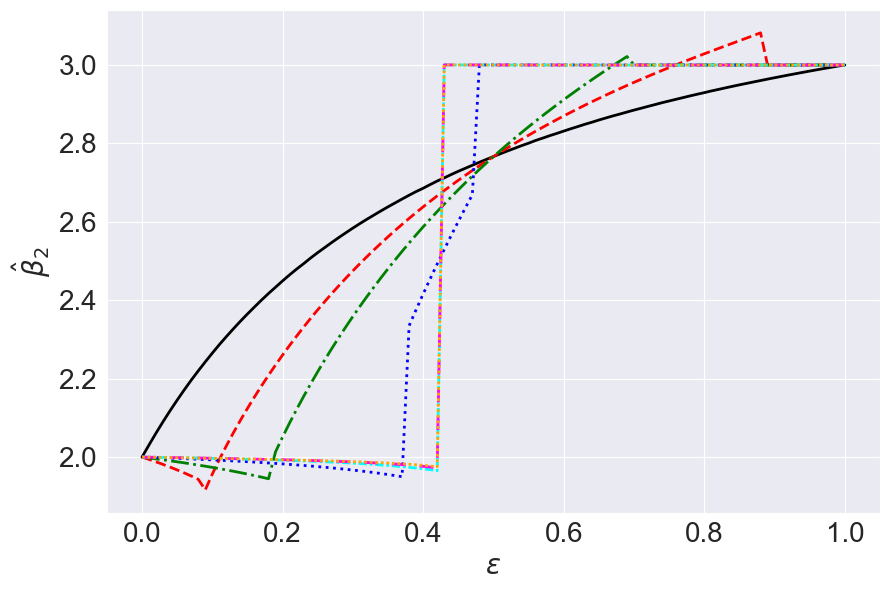}
         \caption{MDPDF of $\beta_2$}
         \label{fig:bp-mm-b2}
     \end{subfigure}
     \hfill
     \begin{subfigure}[b]{0.45\textwidth}
         \centering
         \includegraphics[width=\textwidth]{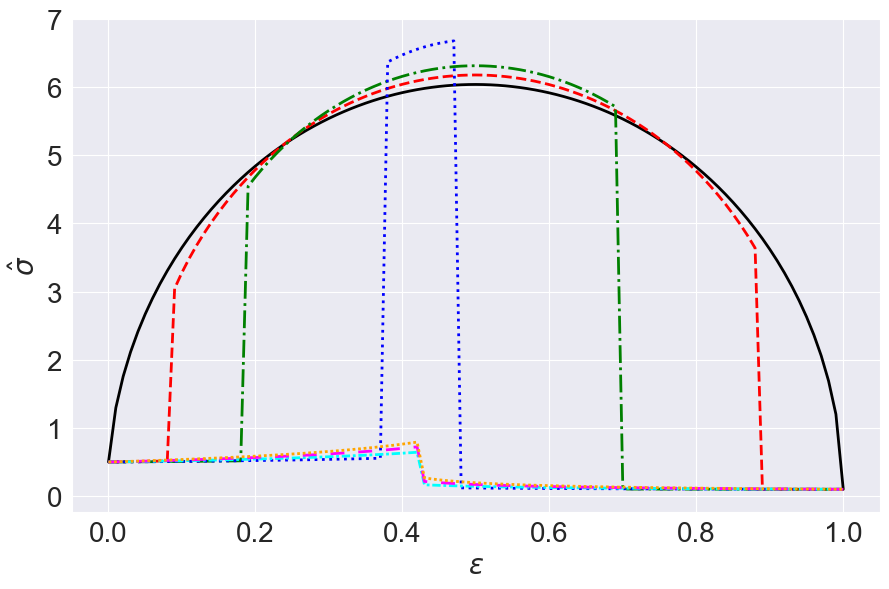}
         \caption{MDPDF of $\sigma$}
         \label{fig:bp-mm-sig}
     \end{subfigure}
     \hfill
     \begin{subfigure}[b]{0.45\textwidth}
         \centering
         \includegraphics[width=0.4\textwidth]{index.png}
         \caption*{}
     \end{subfigure}
     \caption{The MDPDFs obtained under the MM model against various contamination proportions (with a fixed contamination from the model family)}
     \label{fig:bp-mm}
\end{figure}

The visuals in Figure \ref{fig:bp-mm} lead to similar interpretations as in Figure \ref{fig:bp-lr}. At $\alpha=0$, the BPs of the MDPDFs are observed to be zero, and the BPs increase with increasing values of $\alpha$. The MDPDFs of all three parameters are observed to break down near $\epsilon=0.10, 0.19$, and $0.38$ for $\alpha=0.05,~0.10$, and $0.35$, respectively. Further, for all $\alpha\geq 0.5$, the BPs are observed near $0.44$. Here also, the MDPDFs of $\sigma$ create a half-circular pattern for lower values of $\alpha$. \qed
\end{example}

\subsection{Poisson regression} \label{PReg}
Poisson regression is a widely used technique for modeling count data, which models a non-negative integer-valued response with Poisson distributions and relates their mean to the explanatory variables through a log-link function. Let $(y_i,\bm{x}_i),\ i=1,\hdots,n$ be $n$ sample observations on a count response variable $Y$ and covariate vector $\bm{X}=(X_1,\hdots,X_p)'$, where $y_1,\hdots, y_n$ are independent and $x_{ij}$ ($j$-th element of $\bm{x}_i$) are bounded over $i$, for each $j=1,\hdots,p$. Then, under the setup of Poisson regression, the model distribution of $y_i$ is assumed to be Poisson with mean $p_i$ (henceforth denoted as Poisson($p_i$)) and $\ln(p_i) = \bm{x}_i'\bm{\theta}$, where $\bm{\theta}\in \bm{\Theta} = \mathbb{R}^p$ is the model parameter, that needs to be estimated from the data. This is a special case of INH setup with $f_{i,\bm{\theta}}$ denoting the PMF of Poisson($p_i(\bm{\theta})$), with $p_i(\bm{\theta}) = e^{\bm{x}_i'\bm{\theta}}$ for each $i\geq 1$. The MDPDE under this Poisson regression setup is studied in \cite{ghosh2016robust}, along with its asymptotic properties. 

In order to study the ABP of the corresponding MDPDF, let us assume that the true density of $y_i$, denoted by $g_i$, belongs to the parametric model family with parameter value $\bm{\theta}_0$, i.e., $g_i\equiv f_{i,\bm{\theta}_0}$, for all $i$, and $\bm{\theta}_0 \in \bm{\Theta}\setminus \partial\bm{\Theta}$. We also consider the contaminating densities $k_{i,m}$ to be Poisson($\eta_{i,m}$), where $\eta_{i,m}\to\infty$ as $m\to\infty$ for each $i\geq 1$, and $\eta_{i,m}$ satisfies Condition \ref{itm:c2}. Note that, in this case, $\partial\bm{\Theta} = \partial\mathbb{R}^p$ is the set of the points of the form $\bm{\theta}_\infty = ({\theta_{1\infty}, \hdots, \theta_{p\infty}})'$ such that at least one of its components is $+\infty$ or $-\infty$ \sr{(considering boundary with respect to the extended real number system)}. In this case, we consider the breakdown occurs only when any sequence of parameters $\bm{\theta}_m$ (as taken in Conditions \ref{itm:bp2} and \ref{itm:bp4}) $\to\bm{\theta}_\infty$ such that $\bm{x}_i'\bm{\theta}_m \to \infty$, as $m\to\infty$, for all $i$. So, the breakdown is considered only on the subset $\{\bm{\theta}_\infty\in \partial \bm{\Theta}: \bm{\theta}_m \to \bm{\theta}_\infty \mbox{ with }\bm{x}_i'\bm{\theta}_m\to\infty, \mbox{ as } m\to\infty, \mbox{ for all } i\}$. \sr{Furthermore, for the sequence of parameters $\bm{\theta}_m$, we also assume that $p_i(\bm{\theta_m})$ satisfies Condition \ref{itm:c2}.}
% \{\bm{\theta}_\infty: ~\mbox{at least one component of } \bm{\theta}_\infty \mbox{ is } +\infty \mbox{ or } -\infty \mbox{ and the rests belong to the extended real number system}

In this case, Assumptions \ref{itm:bp1}-\ref{itm:bp3} and \ref{itm:bp6} are satisfied with $C = 1$, as shown in Appendix \ref{VAEReg}. So, by Corollary \ref{bp-th-3}, a lower bound of the ABP can be obtained as $\min\{\epsilon_\alpha^\prime, 1/2\}$, where $\epsilon_\alpha^\prime$ is the solution of the equation
\begin{equation} \label{bpeqnpreg}
    \frac{1}{\alpha}x^{1+\alpha} + q_\alpha(1-x) L_0 = 0.
\end{equation}
Now, the MDPDFs of the Poisson regression model parameters at the contaminated distributions can be obtained \sr{by solving the optimization problem} in \eqref{MDPDF} with
% \vspace{-1em}
\begin{equation} \label{MDPDFPReg}
    H_{n,\alpha}^*(\bm\theta) = \frac{1}{n}\sum_{i=1}^n \left[\sum_{y=0}^\infty f_{i,\bm{\theta}}^{\alpha+1}(y) - \left(1+\frac{1}{\alpha}\right) \sum_{y=0}^\infty f_{i,\bm{\theta}}^\alpha(y) g_{i,\epsilon,m}(y) + \frac{1}{\alpha} \right].
\end{equation}
\sr{where $g_{i,\epsilon,m} = (1-\epsilon) f_{i,\bm{\theta}_0} + \epsilon k_{i,m}$ is the $\epsilon$-contaminated density for the $i$-th observation $y_i$}. 

\begin{figure}[H]
     \centering
     \begin{subfigure}[b]{0.45\textwidth}
         \centering
         \includegraphics[width=\textwidth]{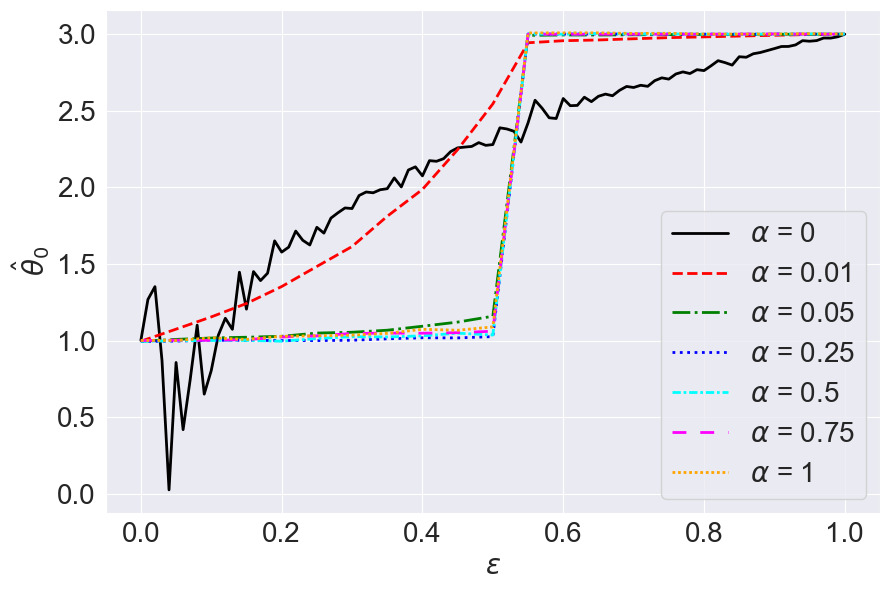}
         \caption{MDPDF of $\theta_0$}
         \label{fig:bp-preg-th0}
     \end{subfigure}
     \hfill
     \begin{subfigure}[b]{0.45\textwidth}
         \centering
         \includegraphics[width=\textwidth]{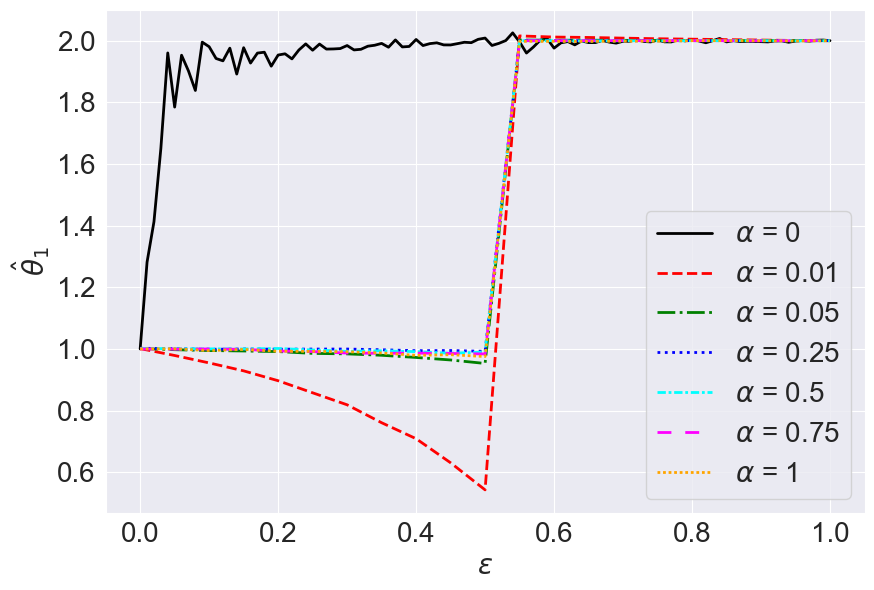}
         \caption{MDPDF of $\theta_1$}
         \label{fig:bp-preg-th1}
     \end{subfigure}
     \caption{The MDPDFs obtained under the Poisson regression model against various contamination proportions (with a fixed contamination from the model family)}
     \label{fig:bp-preg}
\end{figure}

\sr{For numerical illustration, we consider a particular Poisson regression model} given by
\begin{equation} \label{prmodel}
    \ln{p_i} = \theta_0 + \theta_1 x_i\ i=1,\hdots,n,
\end{equation}
and a fixed contaminating distributions as Poisson($3+2x_i$), for each $i\geq 1$. We took the true value of the parameter as $\bm{\theta}_0=(1,1)$, and $x_i$'s were taken as $n=50$ randomly generated observations from Uniform($0,4$) distribution. Under this model setting, for different values of $\alpha$, the MDPDFs of $\theta_0$ and $\theta_1$ were computed at various contamination proportions $\epsilon$, and the resulting estimates are presented in Figure \ref{fig:bp-preg}. To numerically compute the MDPDFs, the inner sums in \eqref{MDPDFPReg} were approximated based on the Monte-Carlo \sr{integration technique by simulating random variables from the underlying Poisson distributions.} Figure \ref{fig:bp-preg} demonstrates that \sr{at smaller values of $\alpha$, namely for $\alpha \leq 0.01$, the MDPDFs are adversely affected by the presence of contamination, and the BPs are observed near $\epsilon = 0$. However, for all $\alpha\geq 0.05$, the MDPDFs are seen to have the highest possible value of BP near $\epsilon = 0.5$.}

\begin{figure}[H]
     \centering
     \begin{subfigure}[b]{0.45\textwidth}
         \centering
         \includegraphics[width=\textwidth]{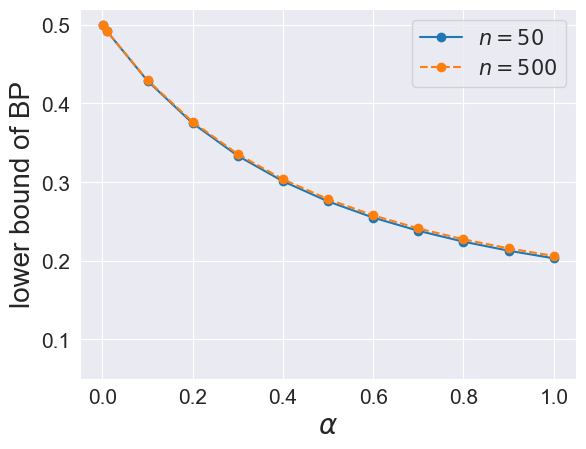}
         \caption{Using Monte-Carlo approximation of $L_0$}
         \label{fig:bp-sol-mc}
     \end{subfigure}
     \hfill
     \begin{subfigure}[b]{0.45\textwidth}
         \centering
         \includegraphics[width=\textwidth]{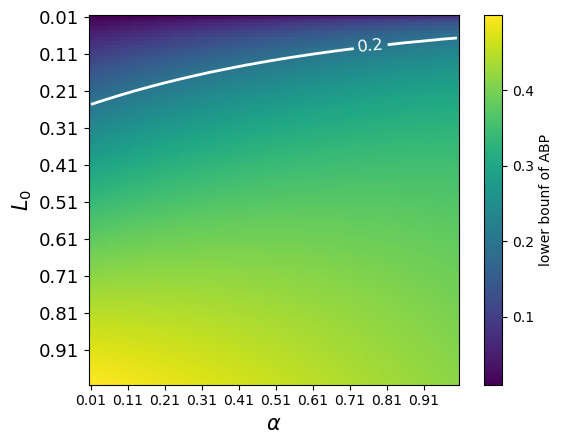}
         \caption{For various values of $L_0$}
         \label{fig:bp-sol-grid}
     \end{subfigure}
     \caption{The lower bound of the ABP, for the Poisson regression setup, obtained by solving \eqref{bpeqnpreg}}
     \label{fig:bp-sol}
\end{figure}

Under the above-mentioned Poisson regression setup, we also computed the lower bound of the ABP of the MDPDFs, as specified by Corollary \ref{bp-th-3}, by numerically solving \eqref{bpeqnpreg}, for different values of $\alpha$ from $0.001$ to $1$. For $n=50$ and $500$, the plots of the lower bounds against various values of $\alpha$ are given in Figure \ref{fig:bp-sol-mc} (where $L_0$ is approximated using the Monte-Carlo integration technique). It shows that the lower bound decreases with increase in $\alpha$ (and is almost independent of the sample size). Yet, the minimum lower bound is $0.20$ at $\alpha=1$, which ensures protection against at least $20\%$ contamination, for any $\alpha$ in $(0,1]$. In Figure \ref{fig:bp-sol-grid}, we have plotted the solutions of Eq.~\eqref{bpeqnpreg} for various fixed values of $L_0$ and $\alpha>0$. It shows that for smaller values of $L_0$, the lower bound of the ABP increases with $\alpha$. We also highlighted the implicit curve with the value of the solution $0.2$, which shows that protection against at least 20\% contamination is ensured for most of the reasonable values of $L_0$ and $\alpha>0$. However, we should note here that the lower bounds presented in the plots may not be very sharp, as $C=1$ is not a sharp upper bound to the inequality in \ref{itm:bp6}. A sharper upper bound $C$ will likely provide more precise lower bounds of the BP, but we do not have any such proof at this moment. We intend to study this problem further in our future research. 

\subsection{Exponential regression} \label{EReg}
This regression setup is a particular case of the Gamma regression, discussed in \cite{cepeda2016gamma}, for modeling non-negative valued response variables (such as lifetime). Suppose there are $n$ sample observations $(y_i,\bm{x}_i),\ i=1,\hdots,n$, on a positive response variable $Y$ and covariate vector $\bm{X}=(X_1,\hdots,X_p)'$, where $y_1,\hdots, y_n$ are independent and $x_{ij}$ ($j$-th element of $\bm{x}_i$) are bounded over $i$, for each $j=1,\hdots,p$. Then, under an exponential regression model (ERM), $y_i$ is modeled by an exponential distribution with mean $p_i$ (henceforth denoted as Exp($p_i$)) and a log-link function $\ln(p_i)=\bm{x}_i'\bm{\theta}$, where the model parameter vector $\bm{\theta}\in \bm{\Theta} (= \mathbb{R}^p)$. Let $f_{i,\bm{\theta}}$ denote the density of Exp($p_i(\bm{\theta})$), with $p_i(\bm{\theta}) = e^{\bm{x}_i'\bm{\theta}}$ for each $i\geq 1$.

Using the definition of the MDPDE for INH data in Section \ref{MDPDE-INH}, for $\alpha>0$, the MDPDE of $\bm{\theta}$ is defined for this setup as $\widehat{\bm{\theta}}_\alpha = \argmin_{\bm{\theta}}H_{\alpha,n} (\bm{\theta})$, where
\begin{equation*}
    H_{\alpha,n}(\bm{\theta}) = \frac{1}{n}\sum_{i=1}^n \left[\frac{1}{(\alpha+1)p_i^\alpha(\bm{\theta})} - \left(1+\frac{1}{\alpha}\right)\frac{1}{p_i^\alpha(\bm{\theta})} e^{-\frac{\alpha y_i}{p_i(\bm{\theta})}} + \frac{1}{\alpha}\right].
\end{equation*}
We now assume that the true density of $y_i$, denoted by $g_i$, belongs to the parametric model family with parameter $\bm{\theta}_0$, i.e., $g_i\equiv f_{i,\bm{\theta}_0}$, for all $i$, and $\bm{\theta}_0 \in \bm{\Theta}\setminus \partial\bm{\Theta}$. We also take the contaminating distributions $K_{i,m}$ as Exp($\eta_{i,m}$), where $\eta_{i,m}\to\infty$ as $m\to\infty$ for all $i$, and $\eta_{i,m}$ satisfies Condition \ref{itm:c2}. Further, note that, for any sequence of parameter $\bm{\theta}_m$, $M_{f_{i,\bm{\theta}_m}} = 1/\{(\alpha+1) e^{\alpha\bm{x}_i'\bm{\theta}_m}\}$, which does not remain bounded if $\bm{x}_i'\bm{\theta}_m \to -\infty$, as $m\to\infty$. So, to avoid this situation, we consider \sr{only exploding-type breakdown, i.e., the sequence} of parameters $\bm{\theta}_m$ (as taken in Conditions \ref{itm:bp2}, \ref{itm:bp4}) satisfy $\bm{x}_i'\bm{\theta}_m \to \infty$, as $m\to\infty$, for all $i$. Here also, the breakdown is considered only in a subset of $\partial\bm{\Theta}$, as mentioned in Section \ref{PReg}. Then, if $p_i(\bm{\theta_m})$ satisfies Condition \ref{itm:c2}, then Assumptions \ref{itm:bp1}-\ref{itm:bp3} and \ref{itm:bp6} are satisfied with $C = 0$, as shown in Appendix \ref{VAEReg}. So, it is evident from Corollary \ref{bp-th-3} and Remark \ref{remark-bp-th-3} that the ABP of the MDPDF of $\bm{\theta}$ is at least $1/2$, for $0<\alpha\leq 1$.

Following Definition \eqref{MDPDF}, the objective function of the MDPDF, at the contaminated distributions with densities $g_{i,\epsilon,m} = (1-\epsilon) f_{i,\bm{\theta}_0} + \epsilon k_{i,m},\ i=1,\hdots,n$, for $\alpha>0$, is given by
\begin{equation*}
    H_{n,\alpha}^*(\bm{\theta}) = \frac{1}{n}\sum_{i=1}^n \left[\frac{1}{(\alpha+1)p_i^\alpha(\bm{\theta})} - \left(1+\frac{1}{\alpha}\right) \left\{\frac{1-\epsilon}{\alpha p_i(\bm{\theta}_0) + p_i(\bm{\theta})} + \frac{\epsilon}{\alpha \eta_{i,m} + p_i(\bm{\theta})}\right\} + \frac{1}{\alpha} \right].
\end{equation*}

\begin{figure}[H]
     \centering
     \begin{subfigure}[b]{0.45\textwidth}
         \centering
         \includegraphics[width=\textwidth]{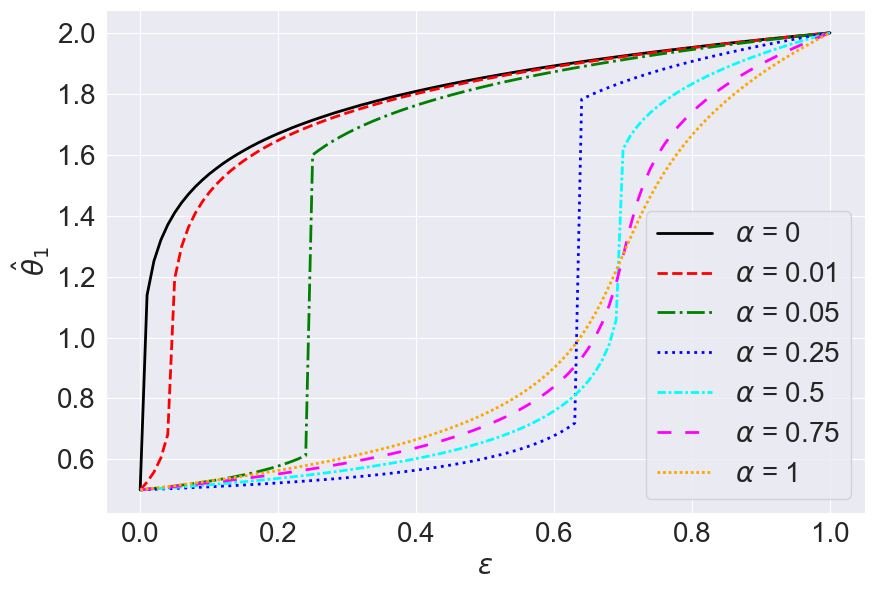}
         \caption{MDPDF of $\theta_1$}
         \label{fig:bp-ereg-th0}
     \end{subfigure}
     \hfill
     \begin{subfigure}[b]{0.45\textwidth}
         \centering
         \includegraphics[width=\textwidth]{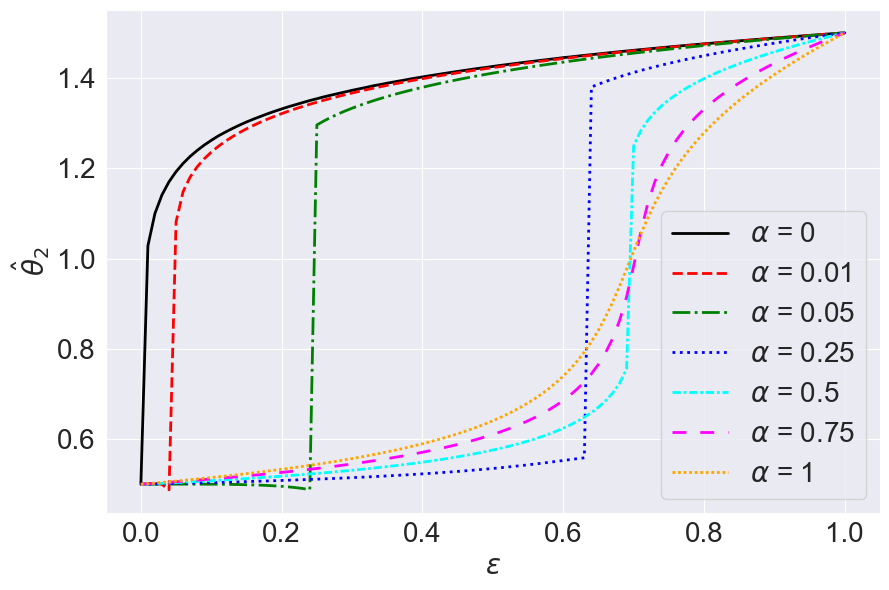}
         \caption{MDPDF of $\theta_2$}
         \label{fig:bp-ereg-th1}
     \end{subfigure}
     \caption{The MDPDFs obtained under the Exponential regression model against various contamination proportions (with a fixed contamination from the model family)}
     \label{fig:bp-ereg}
\end{figure}

For numerical illustration, a particular ERM was considered with the log-link function, which is
\begin{equation} \label{ermodel}
    \ln{p_i} = \theta_1 x_{i1} + \theta_2 x_{i2},\ i=1,\hdots,n,
\end{equation}
and the fixed contaminating distributions is taken to be  Exp($2x_{i1} + 1.5 x_{i2}$), for each $i\geq 1$. The true parameter values were taken as $\bm{\theta}_0=(1/2,1/2)$, and $n=200$ values of the covariates $X_1$ and $X_2$ were generated from $\mathcal{N}(3,1)$ and Uniform($0,5$) distributions, respectively. Under these settings, for different values of $\alpha$, the MDPDFs of $\theta_1$ and $\theta_2$ were computed at various contamination proportions $\epsilon$, and are plotted in Figure \ref{fig:bp-ereg} for different $\alpha\geq 0$. These graphs show that the MDPDFs at $\alpha=0$ and $0.01$ break down as soon as the contamination is added in the underlying distributions, even in a small proportion. At $\alpha=0.05$, the BPs are observed near $\epsilon=0.2$. However, for all choices of $\alpha\geq 0.25$, the MDPDFs are seen to have very high BPs, near 0.5-0.6. 

\section{Simulation studies: Breakdown of the MDPDEs} \label{simulation}
While in Section \ref{Application} we studied the breakdown behavior of the MDPDFs of the parameters of some regression models, in this section we performed some extensive simulation studies to empirically observe the behavior of the (finite-sample) BPs of the corresponding MDPDEs, comparing them with our theoretical results. For each regression model, we utilized the same values of the model parameters, covariates, and contamination scheme mentioned in Section \ref{Application}. The MDPDEs were calculated for $500$ simulated datasets from the underlying regression model at different values of $\alpha\geq 0$ and contamination proportions $\epsilon$. The median values of the resulting MDPDEs were plotted against $\epsilon$, along with the bands corresponding to $25\%$ to $75\%$ quantiles, for different values of $\alpha$ (the band colors were chosen as a transparent version of the colors of the corresponding median MDPDE curves). These plots are provided in Figures \ref{fig:bp-lr-est}-\ref{fig:bp-ereg-est}, for the SLR model \eqref{lr}, MM model \eqref{mmnlr}, Poisson regression model \eqref{prmodel}, ERM \eqref{ermodel}, respectively.  The findings for each case are summarized below.

\begin{itemize}
\item  \textbf{SLR model:} The median MDPDEs of the SLR model parameters (in Figure \ref{fig:bp-lr-est}) show a very similar pattern to those of their respective MDPDFs. The MDPDEs at $\alpha=0$ (i.e., the MLE) have BP near $0$, and it increases as the robustness parameter $\alpha$ is increased. For example, at $\alpha=0.05$ and $0.1$, the BPs are observed in the range of $0.10$-$0.15$ and $0.25$-$0.30$, respectively. Further, at all $\alpha\geq 0.25$, the BPs are observed near $0.45-0.50$. It is also observed that the middle $50\%$ quantile bands get narrower for increasing values of $\alpha$, indicating greater stability of the corresponding MDPDEs.

\begin{figure}[H]
     \centering
     \begin{subfigure}[b]{0.45\textwidth}
         \centering
         \includegraphics[width=\textwidth]{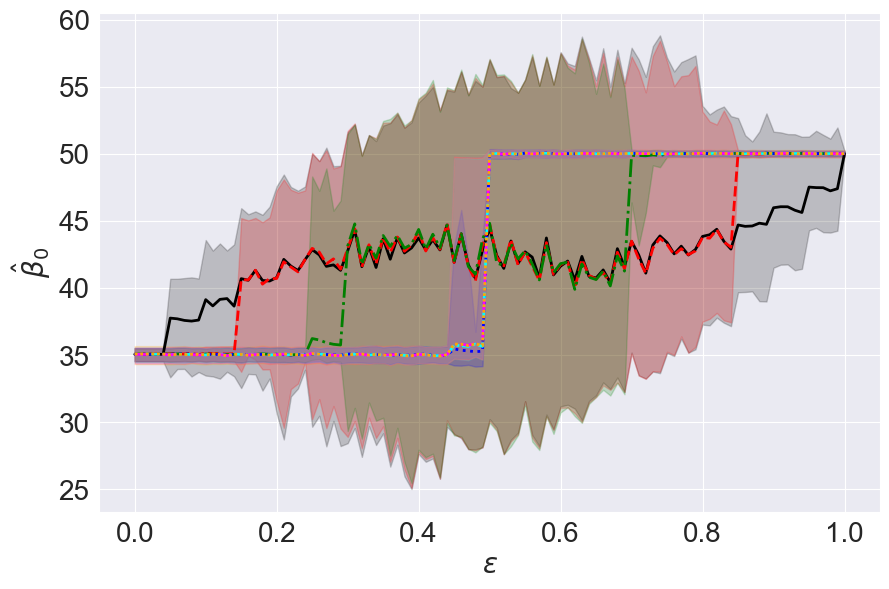}
         \caption{MDPDE of the $\beta_0$}
         \label{fig:bp-lr-b0-est}
     \end{subfigure}
     \hfill
     \begin{subfigure}[b]{0.45\textwidth}
         \centering
          \includegraphics[width=\textwidth]{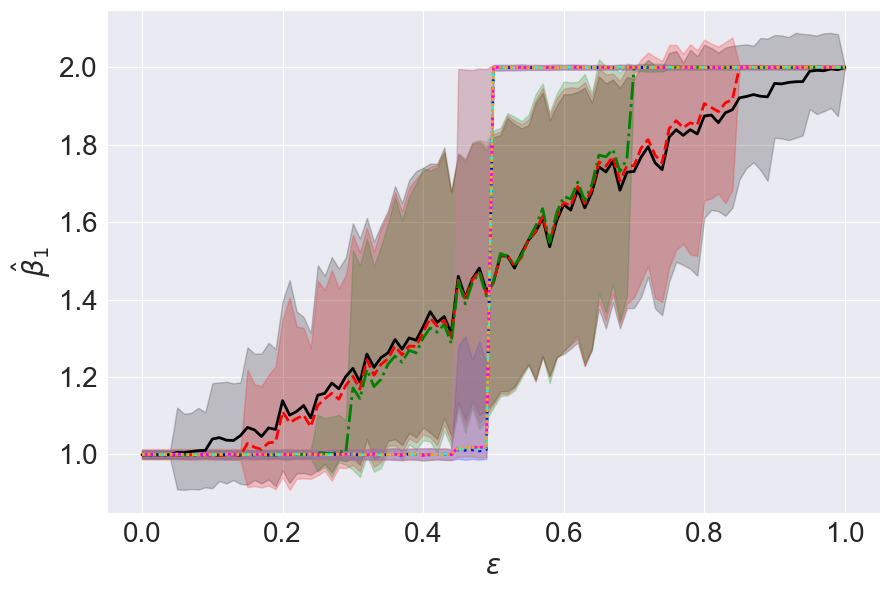}
         \caption{MDPDE of the $\beta_1$}
         \label{fig:bp-lr-b1-est}
     \end{subfigure}
     \hfill
     \begin{subfigure}[b]{0.45\textwidth}
         \centering
         \includegraphics[width=\textwidth]{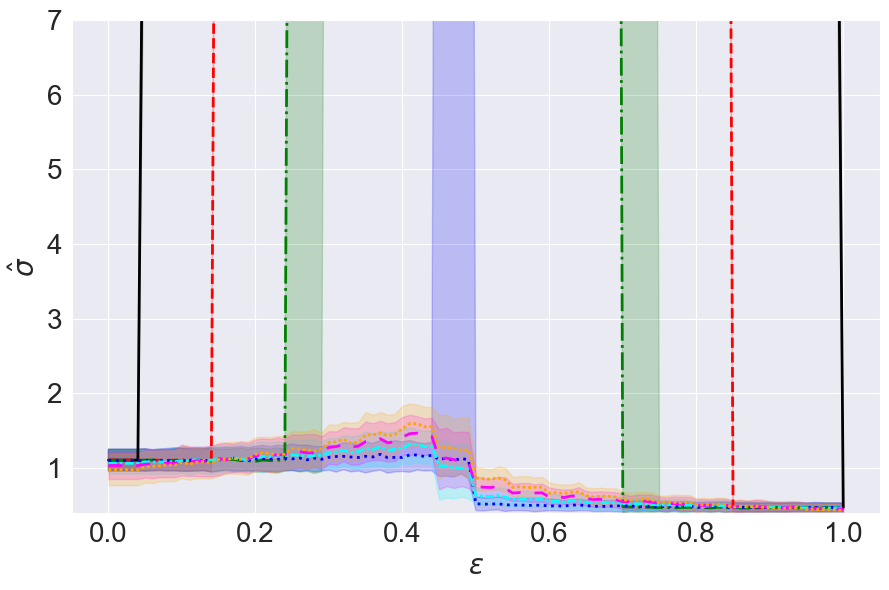}
         \caption{MDPDE of the $\sigma$}
         \label{fig:bp-lr-sig-est}
     \end{subfigure}
     \hfill
     \begin{subfigure}[b]{0.45\textwidth}
         \centering
         \includegraphics[width=0.4\textwidth]{index.png}
         \caption*{}
     \end{subfigure}
     \caption{Plot of median MDPDE and 25\%-75\% quantile band under the SLR model against various contamination proportions}
     \label{fig:bp-lr-est}
\end{figure}

\item \textbf{Non-linear MM model}: The median MDPDEs of the MM model parameters, as plotted in Figure \ref{fig:bp-mm-est}, again behave similarly to those of their respective MDPDFs, under the presence of contamination. At $\alpha=0$, the MDPDEs (i.e., the MLEs) break down as soon as $\epsilon>0$, and the BP increases for increasing $\alpha$. At $\alpha=0.05,~0.10,~0.25$, the BPs are observed near $0.05,~0.15$, and $0.35$, respectively. Further, at all $\alpha\geq 0.5$, the BPs are close to $0.4$. Here also, one can observe that the quantile bands become narrower as $\alpha$ increases, indicating decreasing variability of the corresponding MDPDEs under contamination.

\begin{figure}[H]
     \centering
     \begin{subfigure}[b]{0.45\textwidth}
         \centering
         \includegraphics[width=\textwidth]{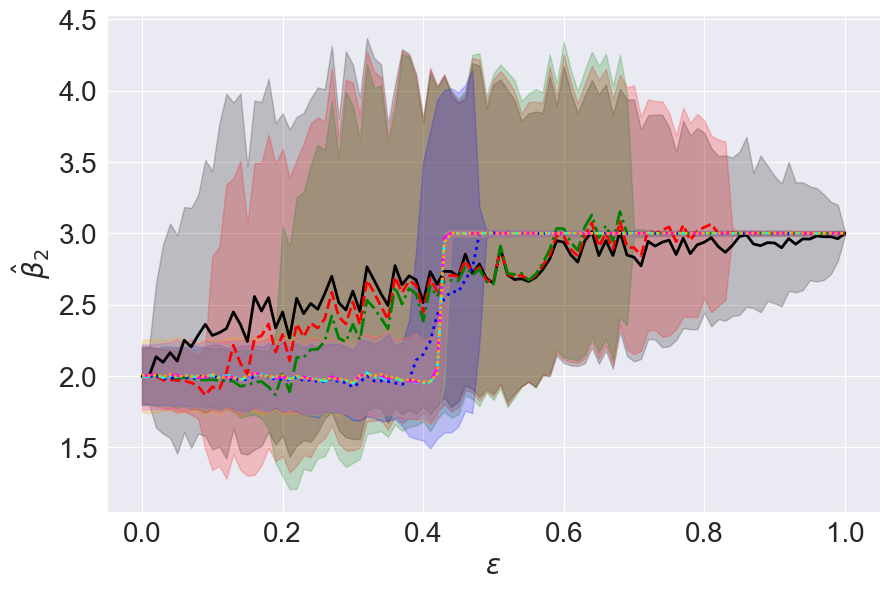}
         \caption{MDPDE of the $\beta_1$}
         \label{fig:bp-mm-b1-est}
     \end{subfigure}
     \hfill
     \begin{subfigure}[b]{0.45\textwidth}
         \centering
         \includegraphics[width=\textwidth]{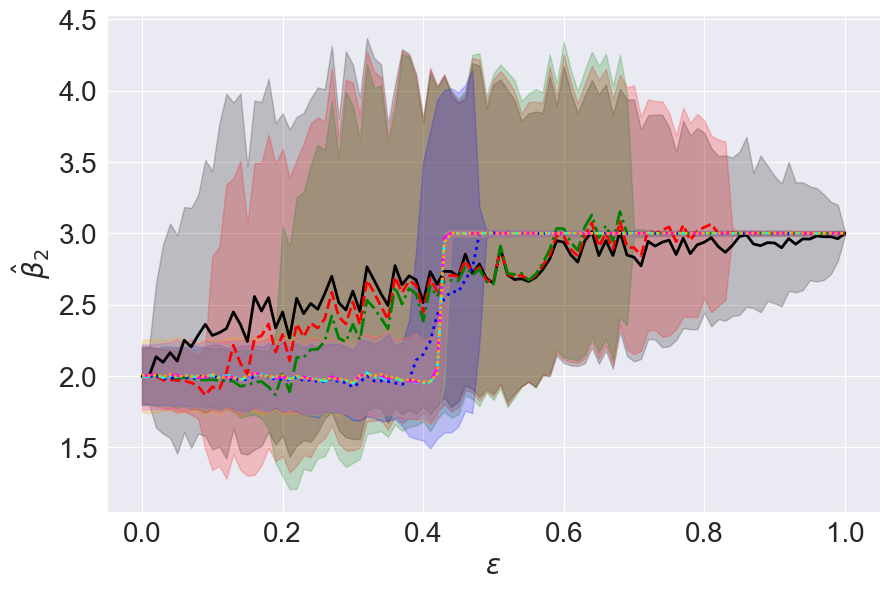}
         \caption{MDPDE of the $\beta_2$}
         \label{fig:bp-mm-b2-est}
     \end{subfigure}
     \hfill
     \begin{subfigure}[b]{0.45\textwidth}
         \centering
         \includegraphics[width=\textwidth]{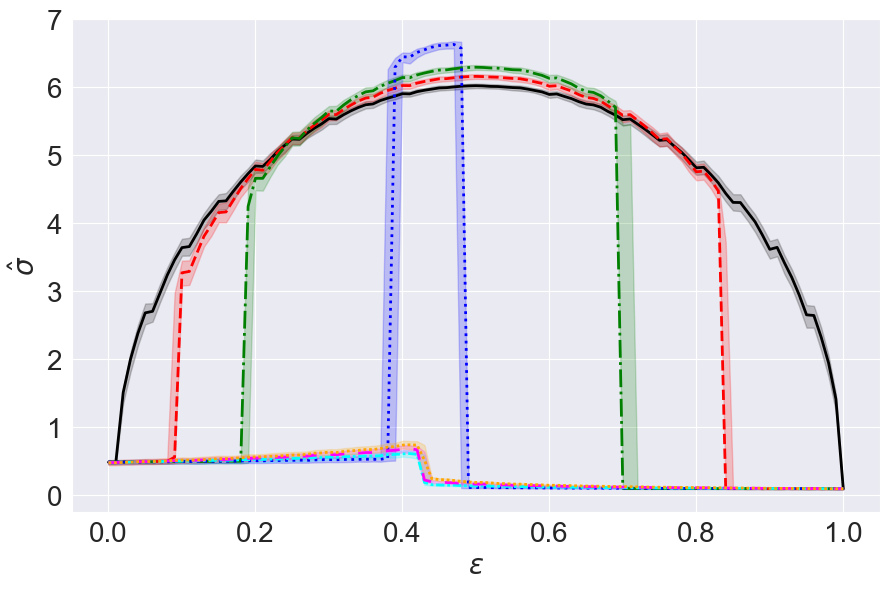}
         \caption{MDPDE of the $\sigma$}
         \label{fig:bp-mm-sig-est}
     \end{subfigure}
     \hfill
     \begin{subfigure}[b]{0.45\textwidth}
         \centering
         \includegraphics[width=0.4\textwidth]{index.png}
         \caption*{}
     \end{subfigure}
     \caption{Plot of median MDPDE and 25\%-75\% quantile band under the MM model against various contamination proportions}
     \label{fig:bp-mm-est}
\end{figure}

\item \textbf{Poisson Regression:} The median MDPDEs of the Poisson regression parameters, as presented in Figure \ref{fig:bp-preg-est}, indicate that higher values of the tuning parameter $\alpha$ prevent the breakdown of the MDPDE until the contamination proportion $\epsilon$ becomes relatively high. At $\alpha=0$ and $0.01$, the BPs are observed near $\epsilon=0$, because these estimates are more adversely affected by the contamination. However, for all $\alpha\geq 0.05$, the BPs are observed near $0.45-0.50$. Additionally, the quantile bands again get narrower as $\alpha$ increases.

\item \textbf{Exponential Regression:} The median MDPDEs of ERM parameters, along with the quantile bands, as presented in Figure \ref{fig:bp-ereg-est}, show a similar pattern to their corresponding MDPDFs, as discussed in the previous section. The MDPDEs at $\alpha=0$ and $0.01$ have the BPs near $\epsilon=0$. At $\alpha=0.05$, the BPs are observed near $\epsilon=0.20$. However, for all $\alpha\geq 0.25$, the MDPDEs are observed to have high BPs (near 0.5-0.6), as in the case for MDPDFs, with the best performance observed at $\alpha=0.25$.

\begin{figure}[H]
     \centering
     \begin{subfigure}[b]{0.45\textwidth}
         \centering
         \includegraphics[width=\textwidth]{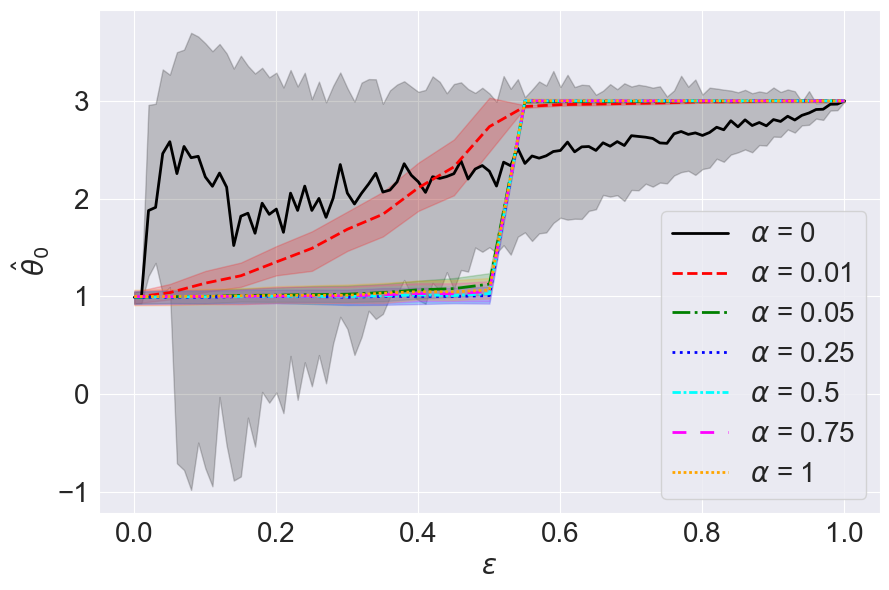}
         \caption{MDPDE of the $\theta_1$}
         \label{fig:bp-preg-th1-est}
     \end{subfigure}
     \hfill
     \begin{subfigure}[b]{0.45\textwidth}
         \centering
         \includegraphics[width=\textwidth]{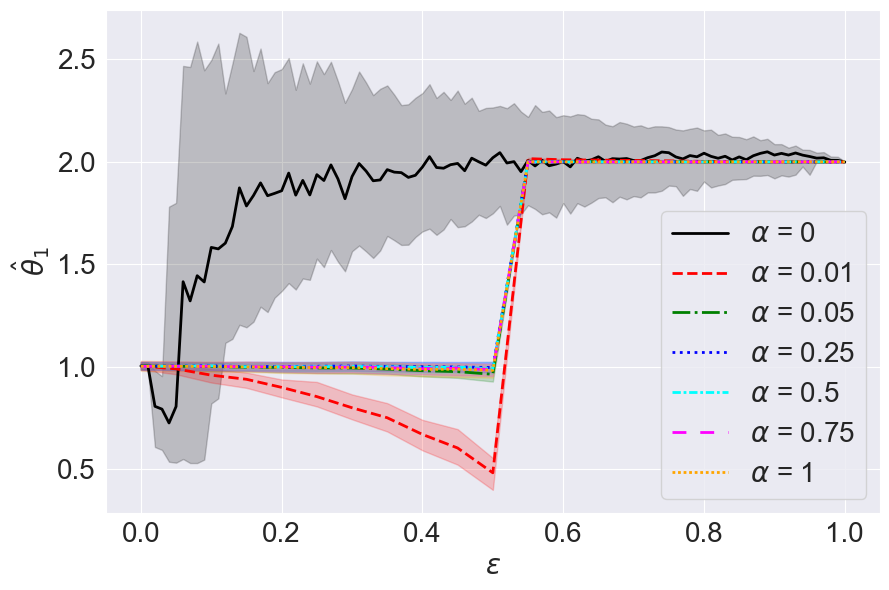}
         \caption{MDPDE of the $\theta_2$}
         \label{fig:bp-preg-th2-est}
     \end{subfigure}
     \caption{Plot of median MDPDE and 25\%-75\% quantile band under the Poisson regression model against various contamination proportions}
     \label{fig:bp-preg-est}
\end{figure}

\begin{figure}[H]
     \centering
     \begin{subfigure}[b]{0.45\textwidth}
         \centering
         \includegraphics[width=\textwidth]{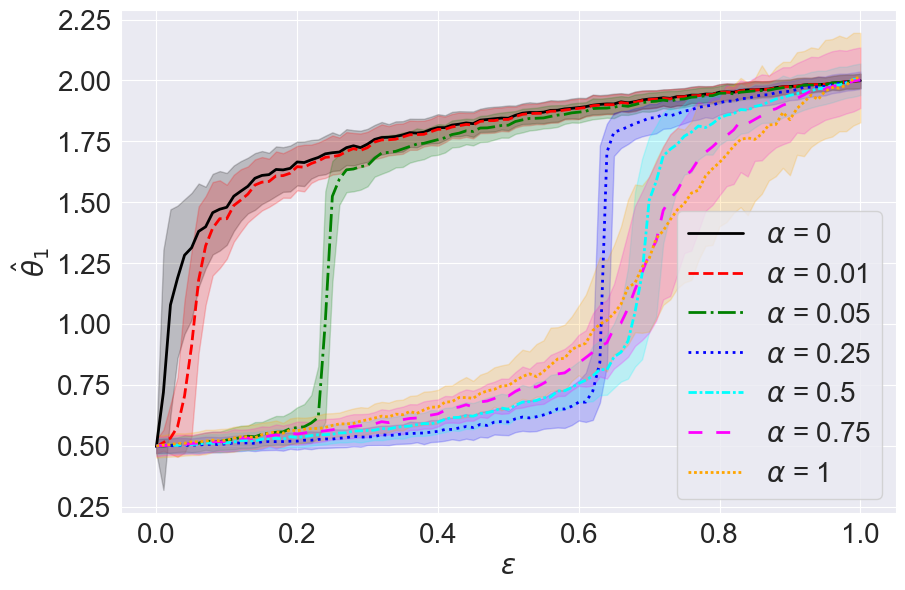}
         \caption{MDPDE of the $\theta_1$}
         \label{fig:bp-ereg-th1-est}
     \end{subfigure}
     \hfill
     \begin{subfigure}[b]{0.45\textwidth}
         \centering
         \includegraphics[width=\textwidth]{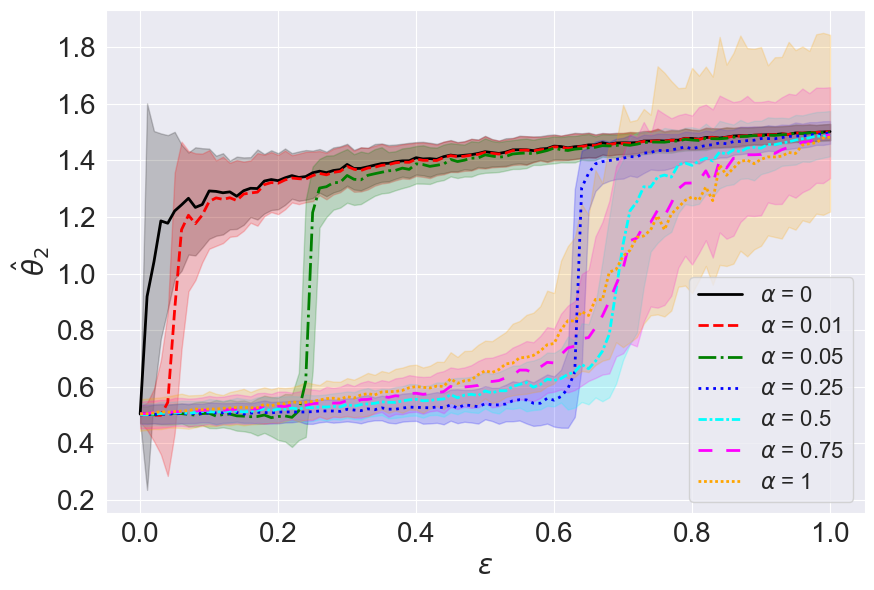}
         \caption{MDPDE of the $\theta_2$}
         \label{fig:bp-ereg-th2-est}
     \end{subfigure}
     \caption{Plot of median MDPDE and 25\%-75\% quantile band under the Exponential regression model against various contamination proportions}
     \label{fig:bp-ereg-est}
\end{figure}
\end{itemize}

\section{Conclusion} \label{conclusion}
In this research, we developed exciting results on the ABP of a well-known and established robust estimator, the MDPDE, under the setup of INH data that constitutes the basic framework for many kinds of (namely fixed-design) regression problems. We provided some precise mathematical definitions of the ABP of an estimator and the associated statistical functional under the INH setup, and subsequently derived a set of conditions leading to the required ABP of the MDPDF. Furthermore, we have established that the ABP of the MDPDE is at least as large as that of the associated functional, although this result relies on some strong assumptions on the nature of the contaminating distributions. However, we believe that this result can also be demonstrated under lighter assumptions, which we intend to explore in future research. Additionally, we conducted extensive simulation studies to observe the behavior of the ABP of the MDPDE and the associated MDPDF, linking them with the general theoretical results. This provides further empirical support and substantiates theoretical development.

This work has further potential to be expanded to more general models and contamination setups, which have, so far, been excluded by our required assumptions. For instance, Assumption \ref{itm:bp2} may not hold when the assumed model distributions of the responses are discrete with the same bounded support (e.g., logistic regression). We hope to develop some such extensions in our future research. It would also be interesting to extend our results to derive the ABP of the MDPDEs under dependent setups, such as time series or stochastic process models.

\appendix
\section*{Appendix}
\section{Proof of the results}\label{Proofs}
\subsection{\textbf{Proof of Theorem} \ref{bp-th}}\label{pf-bp-th}
In this proof, we will use the following Lemma, which was proved in \cite{roy2023asymptotic}.
\begin{lemma} \label{lm1}
     Let $\alpha\geq 0$ and $\{f_m\}_{m\geq 1},\ \{g_m\}_{m\geq 1}$, and $\{\widetilde{g}_m\}_{m\geq 1}$ be three uniformly $L^{1+\alpha}$-integrable sequences of nonnegative functions, i.e., $\sup_m\int f_m^{1+\alpha} d\lambda < \infty$, $\sup_m\int g_m^{1+\alpha} d\lambda < \infty$, and $\sup_m\int \widetilde{g}_m^{1+\alpha} d\lambda < \infty$. Then, if $g_m-\widetilde{g}_m \to 0$ (pointwise) as $m\to \infty$, then we have
    \begin{equation*}
        \left|\int g_m^{1+\alpha} d\lambda - \int \widetilde{g}_m^{1+\alpha} d\lambda\right| \to 0, \mbox{ and } \left|\int f_m^B g_m^A d\lambda - \int f_m^B \widetilde{g}_m^A d\lambda\right| \to 0 \mbox{ as } m\to \infty,
    \end{equation*}
    where $A>0,\ B>0$ and $A + B = 1 + \alpha$.
\end{lemma}

Fix an $\alpha>0$, and define
\begin{equation*}
    D_\alpha(g,f) = f^{1+\alpha} - \left(1+\frac{1}{\alpha}\right)f^\alpha g + \frac{1}{\alpha}g^{1+\alpha},
\end{equation*}
so that $d_\alpha(g,f) = \int D_\alpha(g,f) d\lambda$. Let $\bm{\theta}_m^{(n)} := \bm{T}_\alpha \left(\bm{G}_{\epsilon,m}\right) = \argmin_{\bm{\theta \in \Theta}} \frac{1}{n}\sum_{i=1}^n d_\alpha(g_{i,\epsilon,m}, f_{i,\bm{\theta}})$. Then breakdown occurs for $\bm{T}_\alpha$ at contamination proportion $\epsilon$ if there exist sequences of contaminating densities $\left\{K_{i,m}\right\}_{m\geq 1}$, $i=1,\hdots,n$ and $\bm{\theta}_\infty \in \partial\bm{\Theta}$ such that 
\[\liminf_{n\to \infty} \liminf_{m\to\infty} \lVert\bm{\theta}_m^{(n)}-\bm{\theta}_\infty\rVert = 0,\]
i.e., there exists $N^*\in\mathbb{N}$ such that $\bm{\theta}_m^{(n)} \to \bm{\theta}_\infty$ as $m\to\infty$, for all $n\geq N^*$.

Now, by definition of $\bm{\theta}_m^{(n)}$, we get
\begin{equation}\label{contradict}
    \frac{1}{n}\sum_{i=1}^n d_\alpha \left(g_{i,\epsilon,m}, f_{i,\bm{\theta}_m^{(n)}}\right) \leq \frac{1}{n}\sum_{i=1}^n d_\alpha(g_{i,\epsilon,m}, f_{i,\bm{\theta}_0}), \mbox{ for any } m \mbox{ and } n \geq 1.
\end{equation}
In this proof, we shall find the asymptotic limits of both LHS and RHS of \eqref{contradict} and then deduce a contradiction to imply that there does not exist any sequence of contaminating density $\left\{{k_{i,m}}\right\}$ for which breakdown occurs at $\epsilon \leq \epsilon_\alpha$.\\

\noindent\textbf{Step 1}\\
For a fixed $i\geq 1$, we have,
\begin{equation*}
    d_\alpha\left(g_{i,\epsilon,m}, f_{i,\bm{\theta}_m^{(n)}}\right) = \int_{A_{i,m,n}} D_\alpha\left(g_{i,\epsilon,m}, f_{i,\bm{\theta}_m^{(n)}}\right)d\lambda + \int_{A_{i,m,n}^c} D_\alpha\left(g_{i,\epsilon,m}, f_{i,\bm{\theta}_m^{(n)}}\right)d\lambda,
\end{equation*}
where $A_{i,m,n} = \left\{x: g_i(x) > \max\left(k_{i,m}(x), f_{i,\bm{\theta}_m^{(n)}}(x) \right)\right\}$. Since $g_i = f_{i,\bm{\theta}_0}$ for all $i$, using \ref{itm:bp1} we have, as $m \to \infty$,
\[0\leq\int_{A_{i,m,n}} k_{i,m}d\lambda \leq \int \min\{g_i,k_{i,m}\}d\lambda \to 0, \mbox{ uniformly for each } i, \mbox{ for all } n.\]
Again, an application of \ref{itm:bp2} also gives, as $m\to\infty$,
\[0\leq\int_{A_{i,m,n}} f_{i,\bm{\theta}_m^{(n)}}d\lambda \leq \int \min\{g_i, f_{i,\bm{\theta}_m^{(n)}}\}d\lambda \to 0,\mbox{ uniformly for each } i, \mbox{ for all } n\geq N^*.\]
So, we get, as $m \to \infty$,
\begin{equation*}
    \sup_{1\leq i\leq n} \int_{A_{i,m,n}} k_{i,m}d\lambda \to 0, \mbox{ for all } n, \mbox{ and } \sup_{1\leq i\leq n} \int_{A_{i,m,n}} f_{i,\bm{\theta}_m^{(n)}}d\lambda \to 0, \mbox{ for all } n\geq N^*.
\end{equation*}
Hence, under $k_{i,m}$ and $f_{i,\bm{\theta}_m^{(n)}}$, $A_{i,m,n}$ converges uniformly to a set of probability $0$, as $m \to \infty$, for all $n\geq N^*$. Consequently we have, on $A_{i,m,n}$, $g_{i,\epsilon,m} \to (1-\epsilon)g_i$, as $m\to \infty$, uniformly for each $i=1,\hdots,n$ and for all $n\geq N^*$. Now, note that
\begin{equation*}
    \sup_{1\leq i\leq n} \sup_m \int g_{i,\epsilon,m}^{1+\alpha} \bm{1}_{A_{i,m,n}}d\lambda \leq \sup_{i\in\mathbb{N}} \int g_i^{1+\alpha}d\lambda < \infty, \mbox{ for all } n \mbox{ (by \ref{itm:bp3})}.
\end{equation*}
Then, in Lemma \ref{lm1}, taking $g_m = g_{i,\epsilon,m}\bm{1}_{A_{i,m,n}},\ \widetilde{g}_m = (1-\epsilon)g_i \bm{1}_{A_{i,m,n}}$, and $f_m = f_{i,\bm{\theta}_m^{(n)}} \bm{1}_{A_{i,m,n}}$, (\ref{itm:bp3} ensures the integrability conditions of the lemma) we have
\begin{equation*}
    \sup_{1\leq i\leq n} \left| \int g_{i,\epsilon,m}^{1+\alpha} \bm{1}_{A_{i,m,n}}d\lambda - (1-\epsilon)^{1+\alpha} \int g_i^{1+\alpha} \bm{1}_{A_{i,m,n}}d\lambda\right| \to 0, \mbox{ as } m \to \infty, \mbox{ for all } n \geq N^*,
\end{equation*}
and
\begin{equation*}
    \sup_{1\leq i\leq n} \left| \int f_{i,\bm{\theta}_m^{(n)}}^\alpha g_{i,\epsilon,m} \bm{1}_{A_{i,m,n}}d\lambda - (1-\epsilon) \int f_{i,\bm{\theta}_m^{(n)}}^\alpha g_i\bm{1}_{A_{i,m,n}}d\lambda \right| \to 0, \mbox{ as } m \to \infty, \mbox{ for all } n \geq N^*.
\end{equation*}
Similarly, taking $g_m = f_{i,\bm{\theta}_m^{(n)}} \bm{1}_{A_{i,m,n}},\ \widetilde{g}_m = 0$ and $f_m = g_i$ in Lemma \ref{lm1}, we have
\begin{equation*}
    \sup_{1\leq i\leq n} \int f_{i,\bm{\theta}_m^{(n)}}^{1+\alpha} \bm{1}_{A_{i,m,n}}d\lambda \to 0, \mbox{ and } \sup_{1\leq i\leq n} \int f_{i,\bm{\theta}_m^{(n)}}^\alpha g_i \bm{1}_{A_{i,m,n}}d\lambda \to 0,
\end{equation*}
as $m \to \infty$, for all $n \geq N^*$.

Now, we note that
\begin{align*}
    A_{i,m,n}^c &= \left\{x:g_i(x) \leq k_{i,m}(x)\leq f_{i,\bm{\theta}_m^{(n)}}(x)\right\} \left(=: A_{1,i,m,n}, \mbox{ say}\right)\\
    &\cup\left\{x:g_i(x)\leq f_{i,\bm{\theta}_m^{(n)}}(x) \leq k_{i,m}(x)\right\} \left(=: A_{2,i,m,n}, \mbox{ say}\right)\\
    &\cup\left\{x:k_{i,m}(x)\leq g_i(x) \leq f_{i,\bm{\theta}_m^{(n)}}(x)\right\} \left(=: A_{3,i,m,n}, \mbox{ say}\right)\\
    &\cup\left\{x:f_{i,\bm{\theta}_m^{(n)}}(x) \leq g_i(x) \leq k_{i,m}(x)\right\} \left(=: A_{4,i,m,n}, \mbox{ say}\right)
\end{align*}
Using \ref{itm:bp1}, for $j = 1, 4$, as $m\to\infty$, we have
\[\int g_i \bm{1}_{A_{j,i,m,n}}d\lambda \leq \int \min\{g_i,k_{i,m}\}d\lambda \to 0, \mbox{ uniformly for each } i, \mbox{ for all } n,\]
and using \ref{itm:bp2}, for $j=2,3$, as $m\to\infty$, we have
\[\int g_i\bm{1}_{A_{j,i,m,n}}d\lambda \leq\int\min\{g_i, f_{i,\bm{\theta}_m^{(n)}}\}d\lambda \to 0, \mbox{ uniformly for each } i, \mbox{ for all } n\geq N^*.\]
Hence, we have, as $m\to \infty$,
\begin{equation} \label{bps12}
    0 \leq \int g_i \bm{1}_{A_{i,m,n}^c}d\lambda \leq \sum_{j=1}^4 \int g_i\bm{1}_{A_{j,i,m,n}}d\lambda \to 0, \mbox{ uniformly for each } i, \mbox{ for all } n\geq N^*.
\end{equation}
So, $\int g_i \bm{1}_{A_{i,m,n}}d\lambda \to 1$, as $m\to\infty$, uniformly for each $i=1,\hdots,n$, for all $n>N^*$, i.e., under $g_i$, ${A_{i,m,n}}$ uniformly converges to a set with probability $1$. Thus, we also have
\begin{equation*}
    \sup_{1\leq i\leq n} \left|\int g_i^{1+\alpha} \bm{1}_{A_{i,m,n}}d\lambda - \int g_i^{1+\alpha}d\lambda \right| \to 0, \mbox{ as } m\to\infty, \mbox{ for all } n\geq N^*.
\end{equation*}
Then, combining all these results, we have, uniformly for each $i=1,\hdots,n$,
\begin{equation}\label{bps13}
    \left| \int_{A_{i,m,n}} D_\alpha \left(g_{i,\epsilon,m}, f_{i,\bm{\theta}_m^{(n)}}\right)d\lambda - \frac{(1-\epsilon)^{1+\alpha}}{\alpha} M_{g_i}d\lambda\right| \to 0,\ \mbox{ as } m \to \infty, \mbox{ for all } n \geq N^*
\end{equation}

Again, from the inequality in \eqref{bps12}, we have, $\sup_{1\leq i\leq n} \int_{A_{i,m,n}^c} g_id\lambda \to 0$ as $m \to \infty$, for all $n\geq N^*$, i.e., under $g_i$, $A_{i,m,n}^c$ converges uniformly to a set of probability $0$, as $m \to \infty$. Thus, $\sup_{1\leq i\leq n} \left|g_{i,\epsilon,m} - \epsilon k_{i,m}\right|\bm{1}_{A_{i,m,n}^c} \to 0$, as $m\to\infty$, for all $n\geq N^*$. Also, we have
\begin{align*}
    &\sup_{1\leq i\leq n} \sup_m \int g_{i,\epsilon,m}^{1+\alpha} \bm{1}_{A_{i,m,n}^c}d\lambda\\
    &\leq \sup_{1\leq i\leq n} \sup_m \int \max \left( k_{i,m}^{1+\alpha}, f_{i,\bm{\theta}_m^{(n)}}^{1+\alpha} \right)d\lambda\\
    &\leq \sup_{i\in\mathbb{N}} \sup_m \int k_{i,m}^{1+\alpha}d\lambda + \sup_{1\leq i\leq n} \sup_m \int f_{i,\bm{\theta}_m^{(n)}}^{1+\alpha}d\lambda < \infty, \mbox{ for all } n\geq N^*, \mbox{ due to } \ref{itm:bp3}.
\end{align*}
Then, in Lemma \ref{lm1}, taking $g_m = g_{i,\epsilon,m}\bm{1}_{A_{i,m,n}^c},\ \widetilde{g}_m = \epsilon k_{i,m} \bm{1}_{A_{i,m,n}^c}$, and $f_m = f_{i,\bm{\theta}_m^{(n)}} \bm{1}_{A_{i,m,n}^c}$, we get
\begin{equation*}
    \sup_{1\leq i\leq n} \left| g_{i,\epsilon,m}^{1+\alpha} \bm{1}_{A_{i,m,n}^c}d\lambda - \epsilon^{1+\alpha} \int k_{i,m}^{1+\alpha} \bm{1}_{A_{i,m,n}^c}d\lambda \right| \to 0, \mbox{ as } m \to \infty, \mbox{ for all } n\geq N^*,
\end{equation*}
and
\begin{equation*}
    \sup_{1\leq i\leq n} \left|f_{i,\bm{\theta}_m^{(n)}}^\alpha g_{i,\epsilon,m} \bm{1}_{A_{i,m,n}^c}d\lambda - \epsilon \int f_{i,\bm{\theta}_m^{(n)}}^\alpha k_{i,m} \bm{1}_{A_{i,m,n}^c}d\lambda \right| \to 0, \mbox{ as } m \to \infty, \mbox{ for all } n\geq N^*.
\end{equation*} 
Then combining the above two limits, we have, uniformly for all $i=1,\hdots,n$,
\begin{equation}\label{bps14}
    \left| \int_{A_{i,m,n}^c} D_\alpha\left(g_{i,\epsilon,m}, f_{i,\bm{\theta}_m^{(n)}}\right)d\lambda - \int D_\alpha\left(\epsilon k_{i,m}, f_{i,\bm{\theta}_m^{(n)}}\right)d\lambda \right| \to 0, \mbox{ as } m\to\infty, \mbox{ for all } n\geq N^*.
\end{equation}
Then, combining \eqref{bps13} and \eqref{bps14}, averaging over $i=1,\hdots,n$, for sufficiently large $n,\ (n\geq N^*)$, we have
\begin{align*}
    \liminf_{m\to\infty} \frac{1}{n}\sum_{i=1}^n d_\alpha\left(g_{i,\epsilon,m}, f_{i,\bm{\theta}_m^{(n)}}\right) &= \frac{(1-\epsilon)^{1+\alpha}}{\alpha} \frac{1}{n}\sum_{i=1}^n M_{g_i} + \liminf_{m\to\infty} \frac{1}{n}\sum_{i=1}^n \int D_\alpha\left(\epsilon k_{i,m}, f_{i,\bm{\theta}_m^{(n)}}\right)d\lambda\\
    &= \frac{(1-\epsilon)^{1+\alpha}}{\alpha} \frac{1}{n}\sum_{i=1}^n M_{g_i} + \liminf_{m\to\infty} \frac{1}{n}\sum_{i=1}^n d_\alpha\left(\epsilon k_{i,m}, f_{i,\bm{\theta}_m^{(n)}}\right)\\
    & =: a_{1,n}(\epsilon) \mbox{ (say)}.
\end{align*}

\noindent\textbf{Step 2}\\
Let us again partition the sample space as $B_{i,m} \cup B_{i,m}^c$, where $B_{i,m} := \left\{x: k_{i,m}(x) > f_{i,\bm{\theta}_0}(x) \right\}$. Since $g_i=f_{i,\bm{\theta}_0}$, for all $i\geq 1$, so, using \ref{itm:bp1}, similarly to Step 1, we have
\[\sup_{i\in\mathbb{N}} \int_{B_{i,m}} g_id\lambda \to 0, \mbox{ and } \sup_{i\in\mathbb{N}} \int_{B_{i,m}^c} k_{i,m}d\lambda \to 0, \mbox{ as } m \to \infty.\]
So, under $g_i$ (i.e., $f_{i,\bm{\theta}_0}$), $B_{i,m}$ converges uniformly to a set with probability $0$ and under $k_{i,m}$, $B_{i,m}^c$ converges uniformly to a set with probability $0$. Then, similarly to the previous step, it can be shown that, for all $n$,
\begin{equation*}
    \limsup_{m\to\infty} \frac{1}{n}\sum_{i=1}^n d_\alpha(g_{i,\epsilon,m}, f_{i,\bm{\theta}_0}) = r_\alpha(1-\epsilon) \frac{1}{n}\sum_{i=1}^n M_{g_i} + \frac{\epsilon^{\alpha+1}}{\alpha} \limsup_{m\to \infty} \frac{1}{n}\sum_{i=1}^n M_{k_{i,m}} =: a_{2,n}(\epsilon) (\mbox{say}),
\end{equation*}
where $r_\alpha(\epsilon) = q_\alpha(\epsilon) + \frac{\epsilon^{\alpha+1}}{\alpha}$.\\

Now, note that, due to Assumption \ref{itm:bp4}, $a_{1,n}(\epsilon) - a_{2,n}(\epsilon) > 0$, i.e., $a_{1,n}(\epsilon) > a_{2,n}(\epsilon)$, as $n\to\infty$, which is a contradiction to \eqref{contradict}. So, breakdown cannot occur for any $\epsilon \leq \epsilon_\alpha$. Thus, the asymptotic breakdown point is at least $\epsilon_\alpha$. This completes the proof.

\subsection{Proof of Corollary \ref{bp-th-2}} \label{bp-th-2-pf}
Note that the result is trivially true $\alpha=0$, and so we prove it for any $0<\alpha\leq 1$.

Since, $M_{k_{i,m}}\leq M_{f_{i,\bm{\theta}_m}}$ for all sufficiently large $m$, then by Corollary 3.1 of \cite{roy2023asymptotic}, for all sufficiently large $m$ and for any $\epsilon \leq \frac{\alpha}{1+\alpha}$, we have
\begin{equation*}
    d_\alpha(\epsilon k_{i,m}, f_{i,\bm{\theta}_m})\geq d_\alpha(\epsilon k_{i,m}, k_{i,m}) = \frac{\epsilon^{1+\alpha}}{\alpha} M_{k_{i,m}} + q_\alpha(\epsilon) M_{k_{i,m}}.
\end{equation*}
Now, it is easy to observe that, $q_\alpha(\epsilon) M_{k_{i,m}} - q_\alpha(1-\epsilon) M_{g_i} \geq 0$ as $\epsilon \leq \frac{\alpha}{1+\alpha}$. So, for all sufficiently large $m$ and for any $\epsilon \leq \frac{\alpha}{1+\alpha}$, we have
\begin{equation*}
    d_\alpha(\epsilon k_{i,m}, f_{i,\bm{\theta}_m})\geq \frac{\epsilon^{1+\alpha}}{\alpha} M_{k_{i,m}} + q_\alpha(1-\epsilon) M_{g_i}, \mbox{ for } i=1,2,\hdots.
\end{equation*}
Then averaging over $i=1,2,\hdots,n$ and taking limit as $m\to\infty$ on both sides of the above inequality, \ref{itm:bp4} holds with $\epsilon_\alpha = \frac{\alpha}{1+\alpha}$. Then, the proof is completed by Theorem \ref{bp-th}.

\subsection{Proof of Corollary \ref{bp-th-3}} \label{bp-th-3-pf}
We will prove that the choice of $\epsilon_\alpha^\prime$ given in Corollary \ref{bp-th-3} along with Assumption \ref{itm:bp6} implies Assumption \ref{itm:bp4} with $\epsilon_\alpha = \epsilon_\alpha^\prime$. Note that the LHS of \eqref{bp3ineq} is non-negative. So, it is enough to prove the RHS of \eqref{bp3ineq} is negative.

From Assumption \ref{itm:bp6}, we get, for sufficiently large $n$,
\begin{align*}
    \mbox{ RHS of \eqref{bp3ineq} } &\leq \frac{C\epsilon^{1+\alpha}}{\alpha} + q_\alpha(1-\epsilon) L_0\\
    &= \frac{C\epsilon^{1+\alpha}}{\alpha} + \left\{1-\frac{1+\alpha}{\alpha}(1-\epsilon)\right\} L_0 =: h(\epsilon)\mbox{ (say)}.
\end{align*}
Now, observe that, 
\[h(0) = -\frac{1}{\alpha}L_0 < 0 \mbox{ and } h(1)=\frac{C}{\alpha} + L_0 > 0.\]
Also, $h(\cdot)$ is continuous and strictly increasing in the interval $(0,1)$. So, by Intermediate Value Theorem, the equation $h(\epsilon)=0$ must have a unique root $\epsilon_\alpha^\prime \in (0,1)$ such that $h(\epsilon)<0$ for all $\epsilon<\epsilon_\alpha^\prime$. And hence the RHS of \eqref{bp3ineq} $< 0$ for all $\epsilon<\epsilon_\alpha^\prime$. This completes the proof. 

\subsection{Proof of Corollary \ref{bp-th-4}} \label{bp-th-4-pf}
Since the result holds trivially at $\alpha=0$, so we prove it for any $0<\alpha\leq 1$.

For $\alpha>0$, Assumption \ref{itm:bp4} can be re-specified as
\begin{equation}\label{BP3-RW}
    \frac{1}{n}\sum_{i=1}^n M_{f_{i,\bm{\theta}_m}} - \frac{1+\alpha}{\alpha} \epsilon \frac{1}{n}\sum_{i=1}^n \int f_{i,\bm{\theta}_m}^\alpha k_{i,m}d\lambda > \left\{1-\frac{1+\alpha}{\alpha}(1-\epsilon)\right\} \frac{1}{n}\sum_{i=1}^n M_{g_i},
\end{equation}
for all sufficiently large $m$ and $n$, and for all $\epsilon \in \left[0, \epsilon_\alpha\right)$, where $\epsilon_\alpha \leq \frac{1}{2}$, as given in \ref{itm:bp4}. Note that, RHS of \eqref{BP3-RW} is negative for any $\epsilon<\frac{1}{1+\alpha}$, which is always true as we have $\frac{1}{1+\alpha} > \frac{1}{2}\ (\geq \epsilon)$, for all $0<\alpha\leq 1$. Also, $\frac{1}{n}\sum_{i=1}^n M_{g_i}$ is non-negative and can be arbitrarily large or small, since the true densities $g_i$'s are unknown. So, \eqref{BP3-RW} holds true only when its LHS is positive for sufficiently large $m$, i.e., $\epsilon < \frac{\alpha L}{1+\alpha}$. Thus, \ref{itm:bp4} is satisfied for any $\epsilon < \min\left\{\frac{\alpha L}{1+\alpha}, \frac{1}{2}\right\}$. This completes the proof.

\subsection{\textbf{Proof of Theorem} \ref{bp-th-5}} \label{pf-bp-th-5}
We shall use the contradiction method to prove the theorem. Suppose $\epsilon_2^*<\epsilon_1^*$, and choose an $\epsilon\in (\epsilon_2^*, \epsilon_1^*)$. Then, the MDPDE $\bm{T}_\alpha(\bm{\widehat{G}})$ breaks down at $\epsilon$ proportion of contamination. So, by Definition~\ref{bp-def-E} of the ABP of an estimator, there exists $\bm{\theta}_\infty \in \partial\bm{\Theta}$, $N_1^*\in\mathbb{N}$, and a sequence of contaminating distributions $\{K_m\}_{m\geq 1}$, such that for any fixed $n\geq N_1^*$, $\lVert\bm{T}_\alpha (\bm{\widehat{G}}_{\epsilon,m}) -\bm{\theta}_\infty \rVert \to 0$, as $m\to\infty$, with positive probability, where $\bm{\widehat{G}}_{\epsilon,m} = (\widehat{G}_{1,\epsilon,m}, \hdots, \widehat{G}_{n,\epsilon,m})$ , with $\widehat{G}_{i,\epsilon,m} = (1-\epsilon)\widehat{G}_i + \epsilon K_m,\ i=1,\hdots,n$.

As $\epsilon<\epsilon_1^*$, there will be a contradiction if we can show that
\begin{equation*}
     \liminf_{n\to \infty} \liminf_{m\to\infty} \lVert \bm{T}_\alpha \left(\bm{G}_{\epsilon,m}\right)-\bm{\theta}_\infty \rVert = 0,
\end{equation*}
for some $ \{K_m\}_{m\geq 1}$, where $\bm{G}_{\epsilon,m} = \bm{G}^{(n)}_{\epsilon,m} = (G_{1,\epsilon,m},\hdots, G_{n,\epsilon,m})$ and $G_{i,\epsilon,m} = (1-\epsilon)G_i + \epsilon K_m$. We shall show it in the following two steps.\\

\noindent\textbf{Step 1}

Let us denote $\widehat{\bm{\theta}}(K_m) := \bm{T}_\alpha(\bm{\widehat{G}}_{\epsilon,m})$, and $\bm{\theta}(K_m) := \bm{T}_\alpha(\bm{G}_{\epsilon,m})$. \sr{Invoking Theorem 3.1 of~\cite{ghosh2013robust} under Assumptions (A1)-(A6), we have for any fixed contaminating distribution $K \in \mathcal{K}$,}
\begin{equation*}
    \left\lVert\widehat{\bm{\theta}}(K) - \bm{\theta}(K)\right\rVert \to 0, \text{ as } n\to \infty,
\end{equation*}
\noindent with probability tending to $1$. The above probability convergence is pointwise for every such choice of $K$. Now, consider any two probability distributions $K, K' \in \mathcal{K}$, we get the von-Mises expansion as
\begin{equation*}
    \widehat{\bm{\theta}}(K')
    = \widehat{\bm{\theta}}(K) + \int \text{IF}(x, \widehat{\bm{\theta}}, K)  d(K' - K)(x) + o\left(\Vert K'-K \Vert_{TV} \right),
\end{equation*}
where $\text{IF}(x, \bm{\widehat{\theta}},K)$ is the influence function of the MDPDE $\widehat{\bm{\theta}}$ at the contamination point $x$, when the true distribution is $(1-\epsilon)\bm{G} + \epsilon K$. The last term $o\left( \Vert K'-K \Vert_{TV} \right)$ is the remainder term of the expansion that decays faster than $\Vert K' - K\Vert_{TV}$. Now, an application of Fubini's theorem and the triangle inequality gives
\begin{align}
    \nonumber
    \left\lVert\widehat{\bm{\theta}}(K') - \widehat{\bm{\theta}}(K)\right\rVert &\leq \left[\sup_{x} \left|\text{IF}\left(x, \bm{\widehat{\theta}}, K\right)\right|\right] \int \left| k'(x) - k(x)\right| dx + o\left( \Vert K'-K \Vert_{TV} \right)\\ \label{eqn:von-mises-bound}
    &= 2\left[\sup_{x} \left|\text{IF}\left(x, \bm{\widehat{\theta}}, K\right)\right|\right] \Vert K'-K \Vert_{TV} + o\left( \Vert K'-K \Vert_{TV} \right)
\end{align}
Under the conditions of the theorem, the influence function of the MDPDE is bounded (\sr{c.f. Section 4 of~\cite{ghosh2013robust}}), and hence there exists $M > 0$ such that $\sup_{x} |\text{IF}(x, \bm{\widehat{\theta}}, K)| < M$. As a result, from the inequality \eqref{eqn:von-mises-bound}, we have
\begin{equation*}
    \left\lVert\widehat{\bm{\theta}}(K') - \widehat{\bm{\theta}}(K)\right\rVert \leq O\left( \Vert K' - K\Vert_{TV} \right), ~\mbox{ for all }~ K', K \in \mathcal{K}.
\end{equation*}
\noindent Hence, the family $\left\{\widehat{\bm{\theta}}(K): K \in \mathcal{K}\right\}$ is stochastically equicontinuous. Combining the pointwise consistency with stochastic equicontinuity and using Theorem 2.1 of \cite{newey1991uniform}, we now obtain the uniform consistency of the MDPDE, i.e., with probability tending to $1$,
\begin{equation*}
    \sup_{K \in \mathcal{K}} \left\lVert\widehat{\bm{\theta}}(K) - \bm{\theta}(K)\right\rVert \to 0, \mbox{ as } n\to \infty.
\end{equation*}
\noindent In particular, considering the specific sequence $\{K_m\}_{m\geq 1}$, we have that for any fixed $\epsilon \in [0, 1/2)$ and $\alpha > 0$, $\lVert \bm{T}_\alpha(\bm{\widehat{G}}_{\epsilon,m}) - \bm{T}_\alpha \left(\bm{G}_{\epsilon,m}\right)\rVert \to 0$, as $n\to \infty$, with probability tending to $1$, uniformly for every $m$.\\

\noindent\textbf{Step 2}

Since, the MDPDE $\bm{T}_\alpha(\bm{\widehat{G}})$ breaks down at $\epsilon$ proportion of contamination, for any $\delta>0$, there exists $N_1^* \in \mathbb{N}$ such that $\lVert\bm{T}_\alpha (\bm{\widehat{G}}_{\epsilon,m}) -\bm{\theta}_\infty \rVert < \delta/2$,  whenever $m>m_n$ and $n\geq N_1^*$, with some positive probability $p_0$. Again from Step 1, we obtain $N_2^* \in \mathbb{N}$ such that $\lVert \bm{T}_\alpha (\bm{\widehat{G}}_{\epsilon,m}) -\bm{T}_\alpha (\bm{G}_{\epsilon,m}) \rVert < \delta/2$, whenever $n\geq N_2^*$, with probability at least $1-p_0/2$. So, both of these inequalities hold simultaneously for any $n \geq N^*:=\max\{N_1^*,N_2^*\}$, and any $m>m_n$, with probability at least $p_0/2>0$. So, by triangle inequality, we have
\begin{equation*}
    \left\lVert\bm{T}_\alpha \left(\bm{G}_{\epsilon,m}\right) -\bm{\theta}_\infty \right\rVert \leq \left\lVert \bm{T}_\alpha \left(\bm{\widehat{G}}_{\epsilon,m}\right) -\bm{T}_\alpha \left(\bm{G}_{\epsilon,m}\right) \right\rVert + \left\lVert \bm{T}_\alpha \left(\bm{\widehat{G}}_{\epsilon,m}\right) -\bm{\theta}_\infty \right\rVert \leq \delta,
\end{equation*}
for any $n \geq N^*$, and any $m>m_n$. Since $\delta$ is an arbitrary positive number, the proof follows.

\section{Verification of assumptions for specific regression models}\label{VA}
\subsection{Poisson regression} \label{VAPReg}
Here we shall verify Assumptions \ref{itm:bp1}-\ref{itm:bp3} and \ref{itm:bp6} under the setup of the Poisson regression model, as discussed in Section \ref{PReg}.

At first, it can be easily observed that \ref{itm:bp3} holds as $k_{i,m}, f_{i,\bm{\theta}}$ are probability mass functions of certain Poisson random variables. Now, note that
\begin{align*}
    f_{i,\bm{\theta}}(y) \leq k_{i,m}(y) \Leftrightarrow y \geq \frac{\eta_{i,m}-p_i(\bm{\theta})}{\ln{\eta_{i,m}}-\ln{p_i(\bm{\theta})}} =: v_{i,m}, \mbox{ say}.
\end{align*}
Here, we note that $v_{i,m}\to\infty$ as $m\to\infty$, for each $i$. Also, we have $l_m\leq\eta_{i,m}\leq u_m$, and $p_i(\bm{\theta})$ are bounded in $i$. Let $p^l(\bm{\theta})$ and $p^u(\bm{\theta})$ respectively are some finite lower and upper bound of $p_i(\bm{\theta})$. So, there are sequences $v^l_m = \frac{l_m - p^u(\bm{\theta})}{\ln{u_m}-\ln{p^l(\bm{\theta})}}$ and $v^u_m= \frac{u_m - p^l(\bm{\theta})}{\ln{l_m}-\ln{p^u(\bm{\theta})}}$ such that $v^l_m\leq v_{i,m}\leq v^u_m$, for all $i,m$. So, $v_{i,m}$ is also bounded between two sequences, $v^l_m$ and $v^u_m$, independent of the index $i$. Also, since $l_m\to\infty$, $u_m\to\infty$ as $m\to\infty$ and $\frac{l_m}{u_m} = O(1)$, we have, $v^l_m\to \infty$, $v^u_m\to \infty$ as $m\to\infty$.

Now to check \ref{itm:bp1}, consider the sum
\begin{align*}
    S_{i,m} = \sum_{y=0}^\infty \min \left\{f_{i,\bm{\theta}}(y), k_{i,m}(y)\right\}
    &= \sum_{y=0}^{v_{i,m}} k_{i,m}(y) + \sum_{y = v_{i,m}+1}^\infty f_{i,\bm{\theta}}(y)\\
    &= S_{1,i,m} + S_{2,i,m}, \mbox{ say},
\end{align*}
where $S_{1,i,m},\ S_{2,i,m}$ respectively denote the first and second term of the above sum. Now, we have $S_{2,i,m} \leq \sum_{y=v^l_m}^\infty f_{i,\bm{\theta}}(y) \to 0$ as $m\to\infty$, for all $i$, since $v_m^l\to\infty$ as $m\to\infty$. Also note that,

\begin{align*}
    \sup_i S_{2,i,m} \leq \sum_{y=v^l_m}^\infty \sup_i f_{i,\bm{\theta}}(y) &= \sum_{y=v^l_m}^\infty \sup_i \frac{e^{-p_i(\bm{\theta})}p_i^y(\bm{\theta})}{y!}\\
    &\leq e^{-p^l(\bm{\theta})} \sum_{y=v^l_m}^\infty \frac{(p^u(\bm{\theta}))^y}{y!} \to 0 \mbox{ as } m\to\infty, \mbox{ uniformly for all } i.
\end{align*}
So, $S_{2,i,m} \to 0$ as $m\to\infty$, uniformly for all $i$.

Now, note that
\begin{equation*}
    S_{1,i,m} \leq \sum_{y=0}^{v_{i,m}} k_{i,m}(y) = P\left[X_{i,m} \leq v_{i,m} \right], \mbox{ where } X_{i,m}\sim \mbox{ Poisson}(\eta_{i,m}).
\end{equation*}
Since $\eta_{i,m}\to\infty$ as $m\to\infty$ for all $i$, using central limit theorem (CLT) we have $Y_{i,m}:= \frac{X_{i,m} - \eta_{i,m}}{\sqrt{\eta_{i,m}}} \xrightarrow{D} Z \sim \mathcal{N}(0,1)$ as $m\to\infty$, for all $i$. Now, the moment generating function (MGF) of $Y_{i,m}$ is given by
\begin{equation*}
    M_{Y_{i,m}}(t) = E(e^{tY_{i,m}})
    = e^{-t\sqrt{\eta_{i,m}}} E\left(e^{\frac{tX_{i,m}}{\sqrt{\eta_{i,m}}}}\right) = e^{\frac{t^2}{2}+ \frac{t^3}{3! \sqrt{\eta_{i,m}}}+ \frac{t^4}{4! \eta_{i,m}} + \hdots}
\end{equation*}
Since $l_m\leq\eta_{i,m}\leq u_m$ with $l_m\to\infty$, $u_m\to\infty$ as $m\to\infty$, we have
\begin{equation*}
    \sup_i \left|M_{Y_{i,m}}(t) - M_Z(t)\right| \to 0 \mbox{ as } m \to \infty, \mbox{ for all } t,
\end{equation*}
where $M_Z(t)=e^{\frac{t^2}{2}}$ is the MGF of $\mathcal{N}(0,1)$ distribution. So, using Levy-Cramer CLT (see Chapter 14.7 of \cite{fristedt1996modern}), $Y_{i,m}\xrightarrow{D} Z\sim\mathcal{N}(0,1)$ as $m\to\infty$, uniformly for all $i$.

Now, we have
\begin{equation*}
    P\left[X_{i,m} \leq v_{i,m} \right] = P\left[Y_{i,m} \leq a_{i,m} \right],
\end{equation*}
where
\begin{equation*}
    a_{i,m} = \frac{v_{i,m}-\eta_{i,m}}{\sqrt{\eta_{i,m}}} = \frac{\sqrt{\eta_{i,m}}\left(1-\ln{\eta_{i,m}}+\ln{p_i(\bm{\theta})}\right) - \frac{p_i(\bm\theta)}{\sqrt{\eta_{i,m}}}}{\ln{\eta_{i,m}}-\ln{p_i(\bm{\theta})}}.
\end{equation*}
Note that, $a_{i,m}\to -\infty$ as $m\to\infty$, for all $i$. Also, $a_{i,m}$ is bounded over $i$, for each $m$, since $v_{i,m},\ \eta_{i_m}$ are so. So, there exists sequences $a_m^l$, $a_m^u$ such that $a_m^l \leq a_{i,m}\leq a_m^u$, for all $m,\ i$. It is obvious that $a_m^l\to-\infty$ as $m\to\infty$. Now, we can have
\begin{align*}
    a_m^u = \frac{v_m^u-l_m}{\sqrt{l_m}} = \sqrt{l_m} \left(\frac{\left(\frac{u_m}{l_m}\right)}{\ln{l_m} - \ln{p^u(\bm{\theta})}}-1\right) - \frac{p^l(\bm{\theta})}{\sqrt{l_m}(\ln{l_m} - \ln{p^u(\bm{\theta})})}.
\end{align*}
Since $l_m\to\infty$ as $m\to\infty$ and $u_m/l_m = O(1)$, we have, $a_m^u \to -\infty$, as $m\to\infty$.

Now, for any $\epsilon>0$, there exists $M_\epsilon$ such that $P[Z\leq M_\epsilon]<\epsilon$. Then, there exists some $M^*\in\mathbb{N}$ such that, for all $m>M^*$, we have
\begin{equation*}
    a_{i,m} < M_\epsilon \mbox{ and } \left|P[Y_{i,m}\leq M_\epsilon] - P[Z\leq M_\epsilon]\right| < \epsilon,
\end{equation*}
using the CLT mentioned above. Then, for all $m>M^*$
\begin{align*}
    P[Y_{i,m} \leq a_{i,m}] &\leq P[Y_{i,m} \leq M_\epsilon]\\
    &\leq \left|P[Y_{i,m}\leq M_\epsilon] - P[Z\leq M_\epsilon]\right| + P[Z\leq M_\epsilon]\\
    &\leq 2\epsilon.
\end{align*}
Since $\epsilon$ is arbitrary positive number, we have $P[Y_{i,m} \leq a_{i,m}] \to 0$ as $m\to\infty$, uniformly over all $i$. Hence, $S_{1,i,m} \to 0$ as $m\to\infty$, uniformly over all $i$. So, we have $S_{i,m} \to 0$ as $m\to\infty$, uniformly over all $i$. So, \ref{itm:bp1} is satisfied.

Similarly, \ref{itm:bp2} can be verified by taking $\bm{\theta}=\bm{\theta}_0$ and replacing $\eta_{i,m}$ by $p_i(\bm{\theta}_m) = e^{\bm{x}_i'\bm{\theta}_m}$ in the above derivation. Recall that, $p_i(\bm{\theta}_m)$ also satisfies the same boundedness property \ref{itm:c2} as $\eta_{i,m}$.

It can also be observed that \ref{itm:bp6} holds with $C=1$, for any discrete contaminating distributions.

\subsection{Exponential regression} \label{VAEReg}
Under the Exponential regression setup, as discussed in Section \ref{EReg}, here we check the validity of Assumptions \ref{itm:bp1}-\ref{itm:bp3} and \ref{itm:bp6}.

It is very easy to observe that \ref{itm:bp3} holds. Let us note that,
\begin{equation*}
    f_{i,\bm{\theta}}(y) \leq k_{i,m}(y) \Leftrightarrow y \geq \frac{\ln{\eta_{i,m}}-\ln{p_i(\bm{\theta})}}{\frac{1}{p_i(\bm{\theta})} - \frac{1}{\eta_{i,m}}} =: v_{i,m}, \mbox{ say}.
\end{equation*}
Here we have $v_{i,m}\to\infty$ as $m\to\infty$, for each $i$. Also note that, since $l_m\leq\eta_{i,m}\leq u_m$, and $p_i(\bm{\theta})$ are bounded in $i$, there are sequences $v^l_m\to \infty$, $v^u_m\to \infty$ as $m\to\infty$ such that $v^l_m\leq v_{i,m}\leq v^u_m$, for all $i,m$. So, $v_{i,m}$ is also bounded between two sequences $v^l_m$ and $v^u_m$, independent of the index $i$. Now to check \ref{itm:bp1}, we have
\begin{align}
    \nonumber
    \int_0^\infty \min \left\{f_{i,\bm{\theta}}(y), k_{i,m}(y)\right\}dy &= \int_0^{v_{i,m}} k_{i,m}(y) dy~~~ + \int_{v_{i,m}}^\infty f_{i,\bm{\theta}}(y) dy\\ \label{erb1ch}
    & = 1 - e^{-\frac{v_{i,m}}{\eta_{i,m}}} + e^{-\frac{v_{i,m}}{p_i(\bm{\theta})}}.
\end{align}
Note that, $\frac{v_{i,m}}{\eta_{i,m}} = p_i(\bm{\theta}) \frac{\ln{\eta_{i,m}}-\ln{p_i(\bm{\theta})}}{\eta_{i,m}-p_i(\bm{\theta})} \to 0$ as $m\to\infty$, for each $i$, and $\frac{v_{i,m}}{\eta_{i,m}}$ is bounded over $i$, for each $m$. Since $l_m\leq\eta_{i,m}\leq u_m$, and $p_i(\bm{\theta})$ is also bounded over $i$ (as covariates are bounded), there exists sequences $l^*_m$, $u^*_m$ such that $l^*_m \leq\frac{v_{i,m}}{\eta_{i,m}} \leq u^*_m$, with $l^*_m\to 0$, $u^*_m\to0$ as $m\to\infty$. So, $\frac{v_{i,m}}{\eta_{i,m}} \to 0$ as $m\to\infty$, uniformly for each $i$. So, the integral in \eqref{erb1ch} converges to $0$, as $m\to\infty$, uniformly for each $i$, and uniformly over $\bm{\theta}$ belonging to any compact subset of the parameter space. So, Assumption \ref{itm:bp1} holds.

Similarly, we can show that \ref{itm:bp2} holds by taking $\bm{\theta}=\bm{\theta}_0$ and replacing $\eta_{i,m}$ by $p_i(\bm{\theta}_m) = e^{\bm{x}_i'\bm{\theta}_m}$ in the above derivation, since $p_i(\bm{\theta}_m)$ also satisfies the same boundedness property as $\eta_{i,m}$.

Also, $M_{k_{i,m}} = \int k_{i,m}^{1+\alpha}(y) dy = 1/\{(1+\alpha)\eta_{i,m}^\alpha\} \to 0$ as $m\to\infty$, for each $i$. So, \ref{itm:bp6} holds with $C=0$, for any $\alpha>0$.

\color{black}
\bibliography{Reference.bib}

\end{document}